\documentclass{amsart}
\usepackage{amsaddr}
\usepackage[USenglish]{babel}
\usepackage{amssymb}
\usepackage{graphicx}
\usepackage{amsmath}
\usepackage{color}
\usepackage{enumitem}
\usepackage{latexsym}
\usepackage{hyperref}

\newcommand{\wpp}{{W^{1,p}}}

\newtheorem{theorem}{Theorem}[section]
\newtheorem{lemma}[theorem]{Lemma}
\newtheorem{proposition}[theorem]{Proposition}

\newtheorem{definition}[theorem]{Definition}

\newtheorem{remark}[theorem]{Remark}

\numberwithin{equation}{section}

\title[Fractional Regularity Methods]{On fractional regularity methods for a class of nonlocal problems}

\author[A.L.A. de Araujo]{Anderson L. A. de Araujo}
\address{Departamento de Matem\'atica, Universidade Federal de Vi\c{c}osa, Vi\c{c}osa, MG, 36570-000, Brazil.}
\email{anderson.araujo@ufv.br}

\author[L.H. de Miranda]{Lu\'is H. de Miranda}
\address{Departamento de Matem\'atica, Universidade de Bras\'ilia, Asa Norte, Bras\'\i lia-DF, 70910-900, Brazil.}
\email{demiranda@unb.br}
\thanks{The second author was partially supported by FEMAT-Brazil and FAPDF-Brazil, grants 01/2014 and 4749.25.27523.07072015.}



\begin{document}

\begin{abstract}
In the past years, the phenomenon of fractional regularity has been addressed for a large class of linear and/or quasilinear differential operators, mostly, in terms of certain Besov spaces. As it turned out, for equations governed by the $p$-Laplacian, in general, the regularity of solutions appears in terms of functional spaces with nonlinear order of smoothness. Moreover, despite its own interest, fractional regularity methods may be used as a tool for the investigation of some Partial Differential Equations which are not usually addressed in this manner. Thus, the purpose of the present paper is to exploit such methods in order to provide some results regarding existence and regularity of solutions to a class nonlocal elliptic equations which are linked to the $p$-Laplacian. This is done by means of explicit a priori estimates regarding Lebesgue and Nikolskii spaces, which are part of the present contribution.  As a consequence, this approach allows a relaxation on some of the standard conditions employed in this class of problems.
\end{abstract}

\maketitle


\par
\section{Introduction}
In the present paper, we address some aspects on the existence and fractional regularity of solutions for
\begin{equation}\tag{P}\label{p1}
\begin{cases}
\begin{array}{rll}
	\displaystyle -\left[a\left(\|u\|^p_{1,p}\right)\right]^{p-1}\Delta_p u + u &=& f(x,u) \ \ \mbox{ in } \ \ \Omega\\
\displaystyle	\frac{\partial u}{\partial \eta}&=&0 \ \ \mbox{ on }\ \ \partial \Omega,
	\end{array}
\end{cases}
\end{equation}
where  $\Omega \subset \mathbb{R}^N, \ N \geq 2 ,$  is an open bounded domain with smooth
boundary $ \partial \Omega $,  $ \Delta_p $ is the $p$-Laplacian operator
\[ \Delta_p u  = \text{div}\bigg( |\nabla u |^{p-2} \nabla u\bigg) , \quad \text{ with } p>2 \]
and $a(.)$ is the so--called $p$-Kirchhoff, or Kirchhoff term, which will be assumed to be continuous and bounded by below.
\vspace{0.2cm}

The interest on this sort of nonlocal problem goes back to G. Kirchhoff in the end of the 19$^{th}$ century, cf. \cite{kirchhoff}, and has been addressed by a vast literature after the seminal work due to J.L. Lions in \cite{lions}, where the original hyperbolic version of \eqref{p1} is investigated. Although that it is remarkable that Kirchhoff problems are related to the modelling of the nonlinear vibration phenomenon, for the sake of brevity, we refer the reader to \cite{acg,cog1,cog2,kirchhoff,lions,rak,wxz} and the references thereof for further information on this subject, where the physical background is discussed in more detail.

Moreover, recently, issues concerning the improved regularity of solutions for $p$-Laplacian--like equations and its several generalizations, such as systems or double--phase operators, have also attracted a considerable attention within the field of Partial Differential Equations. Actually, results which guarantee the validity of generalized Calder\'{o}n-Zygmund inequalities, higher order integrability of the gradients, partial regularity results or fractional order regularity for the solutions have been investigated extensively in the late years. Without the intention of being complete, we cite \cite{akg,bcdks,cianchi,coz,duzaar,com1,com2,dmm,dmm2,kum1,kms1,kms2}, and the references thereof, where these issues where addressed by means of a variety of approaches, e.g., from the so--called nonlinear Wolff potentials to direct methods, and so on. One of the byproducts of the contributions listed above relies on the gain of compactness for derivatives of order greater than one.
Most of all, despite that this sort of result is intrinsically interesting, it is indeed the potential for possible applications on a broad class of different contexts of Partial Differential Equations which may reinforce its role.
In particular, we remark that in \cite{coz}, following an approach based on Nirenberg's translation methods, the author provides fractional regularity regarding Nikolskii and Sobolev--Slobodeckii spaces for another sort of nonlocal operators which are close to the fractional $p$-Laplacian. Complementarily, in \cite{kms1,kms2} the authors address fractional regularity results for nonlocal fractional Laplacian--like operators by means of an improved version of the Gehring's inequality of fractional character.

It is our purpose to give an example of such application, i.e., to employ certain fractional order regularity results regarding derivatives of order greater than one to investigate nonlocal problems of Kirchhoff's type, where the integral terms carry out nonlinearities related to gradient terms.

First, by adapting fractional regularity results to the context on Kirchhoff problems, we prove  that solutions to problem \eqref{p1} satisfy a set of a priori bounds, including estimates in the so--called Nikolskii spaces $\mathcal{N}^{1+\frac{2}{r},r}$, a special class of the Besov interpolation spaces, namely $\mathcal{B}^{1+\frac{2}{r}}_{r,\infty}$.
Let us remark that for the reader convenience, we will discuss the relation of these spaces and $p$-Laplacian--like equations later, see Subsection \ref{fractionalspaces} and Section \ref{appendix} below. In short, by taking advantage of the compact embedding
\[\mathcal{N}^{1+\frac{2}{r},r}\hookrightarrow\hookrightarrow W^{1,p}, \mbox{ where }r>2,\]
we are able to conclude that problem \eqref{p1} has a solution under certain basic restrictions to $a(.)$, the nonlocal Kirchhoff coefficient. As a matter of fact, we assume that $a(.)$ is bounded by below and continuous, only. The use of fractional regularity methods to handle the effects of Kirchhoff--like nonlinear terms is new, being the main contribution of the present work. Actually, this is made by means of a priori estimates on Nikolskii spaces provided in \cite{dmm,dmm2} and their adaptations developed in Section \ref{apriori} below.

Further, by using Moser's iteration technique we obtain $L^\infty$ estimates in terms of Sobolev norms of solutions, where
 our approach is based on the assumption of a version of the so--called nonquadraticity conditions for $f$ and its primitive $F(.,.)$. For the best of our knowledge,  this use of nonquadraticity--like conditions is also new and may be applied to other contexts of Partial Differential Equations.

It is worth to notice that one of the main challenges of the present paper is to handle the $p$-linear unbalance between the high and low order terms in the left--hand side of \eqref{p1}, coupled to a Kirchhoff term, which will be assumed to be bounded by below and may possibly fail to be of bounded variation. This brings on several extra difficulties to control the norms of fractional or integer order for the solutions of the equation which is considered in this work.

\textbf{Plan for the paper.}
In Section \ref{notations} we introduce the basic notation, state the main hypotheses and contributions which are going to be discussed throughout the text.
On Section \ref{PreliminaryResults}, we develop certain basic tools addressing convenient energy estimates of fractional order for a linearized version of \eqref{p1},  whereas explicit a priori estimates of the full problem are obtained in Section \ref{apriori}. Further, in Section \ref{proof}, we exhibit the proofs of the main results of this paper and finally, in Section \ref{appendix}, we present the basic background on Fractional Regularity used in the text. This includes an example which illustrates the role of Fractional Regularity methods in the investigation of $p$-Laplacian--like equations.

\section{Notations, general hypotheses and main results}\label{notations}

In this paper,
$\Omega \subset \mathbb{R}^N, \ N \geq 2,$ is considered to be an open bounded  domain  with smooth
boundary $ \partial \Omega ,$ which means that $\partial \Omega \in C^\infty$, and its unit normal vector will be denoted by $\eta$. Moreover, throughout the present discussion, $a.e.$ will always be with respect to $\mu$, the Lebesgue $N$-dimensional measure.

For the sake of clarity of the estimates, inequalities or calculations which will be found throughout this paper, we are going to introduce some basic notation. Indeed, we denote by $p^*$ the standard Sobolev critical exponent:
\[
p^*=\dfrac{Np}{N-p}, p<N, \mbox{ or } +\infty \mbox{ otherwise.}
\]

In addition, we set
\begin{equation}
\label{exponents2}
r_s=s(p-2)+2, \mbox{ where } s \in [2,+\infty) ,
\end{equation}
which is the basic regularity exponent of the present work. This exponent is known to be associated to $p$-Laplacian-like equations and appears naturally in the context of Fractional Regularity theory for such operators, see \cite{dmm,dmm2}. Remark that for the limit case $p=2$, $r_s=2$ for all $s\geq 2$.

Moreover, in order to abbreviate certain exponents found in the a priori estimates discussed below, see Sections \ref{PreliminaryResults} and \ref{apriori}, we introduce the quotient
\begin{equation}
\label{exponents3}
\lambda=\dfrac{p^*}{p}, \mbox{ if } p<N, \mbox{ or } +\infty \mbox{ otherwise},
\end{equation}
and, for $p>2$, $\tau>2$, the function
\begin{equation}
\label{exponents4}
\omega_{p,\tau} =
\begin{cases}
0 &\mbox{ if } p\geq 3-\dfrac{2}{\tau} \\
\dfrac{\big(\tau(3-p)-2\big)(p-1)}{\big(\tau(p-2)+2\big)(p-2)} &\mbox{ if } 2<p<3-\dfrac{2}{\tau}.
\end{cases}
\end{equation}

At this point, we must stress that despite the extremely nonlinear appeal of the last exponents, the later choices make our estimates considerably more concise, what justifies their introduction.

In addition, for the sake of clarity of the a priori or energy estimates presented in the text, we introduce a class of auxiliary polynomials. In fact, by considering
\begin{align}
&\label{4lm2}
n_s=\min\{n\in \mathbb{N}:p^* (k_n+1)\geq s\}, 
\end{align}
for a given nondecreasing nonnegative sequence $\{k_n\}_\mathbb{N}$, we define the polynomials
\begingroup
\Large
\begin{equation}
\label{2lm2}
 h_s(x)= x+x^{(n_s+1)\lambda},x\geq0.
\end{equation}
\endgroup
Despite being rather technical, the introduction of these polynomials will aid us establishing uniform notation for most of our a priori estimates, see Propositions \ref{prop2}, \ref{prop3} and \ref{prop4}, on Section \ref{apriori}.

We also stress that we are going to use the standard notation for embeddings between Banach spaces. Actually, given two Banach spaces $A$ and $B$,  $ A \hookrightarrow B$ means that $A$ continuously embedded in $B$. In addition, by $A \hookrightarrow\hookrightarrow B$, we mean that $A$ is compactly embedded in $B$.
\subsection{Fractional Spaces and the $p$-Laplacian}
\label{fractionalspaces}
Throughout our entire discussion, we make use of spaces of fractional order of differentiability, and some of its particular properties. Indeed, the so--called Nikolskii spaces \[\mathcal{N}^{\sigma,r}, \mbox{ for  } 1<\sigma<2 \mbox{ and }r\geq 2,\]
with norm
\begin{equation}
\nonumber
\|u\|_{\mathcal{N}^{\sigma,r}}=\sum_{i=1}^N\sup_{h\neq 0} \bigg{(}\int_{\Omega_{|h|}}\dfrac{\big{|}\partial_{x_i} u(x+h)-\partial_{x_i} u(x)\big{|}^r}{|h|^{\sigma r}}\bigg{)}^{1/r} +\|u\|_{W^{1,r}},
\end{equation}
will be used along the text. It is remarkable that this sort of space, as well as other kinds of fractional spaces, such as Sobolev--Slobodeckii spaces, could be considered as intermediate spaces between the Sobolev spaces $W^{1,r}$ and $W^{2,r}$.

In spite of the fact that there are several ways to introduce fractional spaces, e.g., by means of interpolation or direct methods,  it is beyond our purpose to describe in detail their specific properties -- we solely intent to employ what is convenient to our approach. Anyway, for the sake of completeness, we have collected part of the essential background on these spaces in order to work with Fractional Regularity theory for PDE's. We although suggest for the interested reader the references \cite{bin,kufner,leoni,leoni2,npv,rsc,saw} or \cite{triebel3}, where such spaces are discussed thoroughly. All in all, for our purposes, the main property regarding the fractional spaces considered along the text is the following compact embedding
\[
\mathcal{N}^{\sigma,r} \hookrightarrow\hookrightarrow W^{1,r}, \mbox{ for  } 1<\sigma<2,
\]
 see Lemma \ref{compactness} in the Appendix. Basically, this result is the key feature allowing us to relax the hypotheses on the Kirchhoff term. 

It is a well--known issue the loss of regularity for the solutions of quasilinear PDE's, with very distinct patterns of regularity between the singular and degenerate cases. In turn, regarding the degenerate one, in general, the solutions to $p$-Laplacian--like equations, for $p>2$, lack regularity for the integer derivatives of order greater than one. Indeed, for data $f\in L^s$, in general, the solutions to this sort of equation do not belong to the spaces $W^{2,r_s}$ or $W^{2,p}$, since for $r_s$ or $p$ sufficiently large, this  would contradict their $C^{1,\alpha}$  threshold regularity, see for instance \cite{dibenedetto2,tolks} or the example in the Appendix, see p. \pageref{exampleapp}. Thus, in order to seek for higher regularity for the solutions of such equations, it is natural to consider spaces which are intermediate between $W^{1,p}$ and $W^{2,p}$ or $W^{2,r_s}$. As it turned out, Nikolskii spaces are suitable to fill this gap, for instance see \cite{com1,coz,ebmeyer,ebmeyer2,garroni,dmm,dmm2,simon1,simon2}. This is an outcome of the intrinsic energetic estimates for solutions of $p$-Laplacian--like equations. Indeed, to illustrate the ideas and to avoid unnecessary additional technical details, let us consider  the case $s=2$ so that $r_s=\frac{2p}{p^\prime}$. We stress that, given a solution to
\begin{equation}
\nonumber
\begin{cases}
\begin{array}{rll}
-\Delta_p u+u &=& f \mbox{ in } \Omega
\\
\dfrac{\partial u}{\partial \eta} & =&0 \mbox{ on } \partial \Omega,
\end{array}
\end{cases}
\end{equation}
it is possible to prove that it necessarily belongs to $\mathcal{N}^{p^\prime, \frac{2p}{p^\prime}}$. In fact, in analogy to the linear case, by considering the most natural test function, $-\Delta_p u$, it is clear that
\[\int_\Omega |\Delta_p u|^2\leq \|f\|^2_{L^2}\]
and then, by integrating by parts twice, after some straightforward manipulations, formally, we arrive at
\[\int_\Omega |\nabla u|^{\frac{2p}{p^\prime}-2} |D^2 u|^2\leq C \int_\Omega |\Delta_p u|^2 + \mbox{ lower order terms}.\]
 Further, by means of certain nonlinear estimates, see Lemma \ref{normalemma} in the Appendix, we obtain

\[
 [[u]]^{2p/p^\prime}_{\mathcal{N}^{p^\prime,\frac{2p}{p^\prime}}}\leq \int_{\Omega} |\nabla u|^{\frac{2p}{p^\prime}-2}|D^2u|^2,
\]
and thus, after a standard interpolation combined with the embeddings given in Lemma \ref{compactness}, we obtain that
\[\|u\|^{2p/p^\prime}_{\mathcal{N}^{p^\prime,\frac{2p}{p^\prime}}}\leq C\big{(}\|f\|^2_{L^2}+1\big{)}.\]

Due to the latter estimates, which, we stress, are intrinsic to equations involving the $p$-Laplacian, we are led to the framework of Nikolskii spaces $\mathcal{N}^{\sigma,r}$, or at least, to another sort of spaces of fractional order which are somehow close to this choice. For instance, let us mention that another possibility would be the standard Sobolev spaces of fractional order,  i.e., the Sobolev--Slobodeckii spaces $W^{\sigma,r}$, which are topologically close to $\mathcal{N}^{\sigma,r}$ in the sense that
\[ \mathcal{N}^{\sigma,r}\hookrightarrow\hookrightarrow W^{\sigma-\epsilon,r} \hookrightarrow \mathcal{N}^{\sigma-\epsilon,r},\]
for $\epsilon>0$ sufficiently small, see Lemma \ref{compactness}. Nevertheless, due to the nature of the estimates which we employ along the text, the choice of  $\mathcal{N}^{\sigma,r}$ turned out to be more convenient.
\subsection{Main Results}
\label{mainresults}

First, motivated by the latter remarks, let us start by stating our basic hypotheses.
In fact, the nonlocal term $a: \mathbb{R} \rightarrow \mathbb{R}$, is supposed to be a continuous function for which there exists $a_0>0$ satisfying
\begin{equation}
a(t)\geq a_0, \mbox{ for every } t\geq 0.
\tag{H1}\label{H1}
\end{equation}
Once again, with the purpose of visually simplifying part of the subsequent computations, we define the following auxiliary functional
\begin{equation}
\nonumber
b(u)=a(\|u\|^p_{W^{1,p}})^{p-1} \mbox{ so that } b(u)\geq a_0^{p-1}, \mbox{ for } u \in W^{1,p}.
\end{equation}

In addition, the nonlinear term $f(x,t)$ denotes a Caratheodory function possessing subcritical growth and a nontrivial character.
Indeed, if $2<p<N$, we suppose that
\begin{equation}
\begin{cases}
\tag{H2}\label{H2}
&\mu\bigg(\{x\in \Omega: f(x,t)\neq t\}\bigg)>0 \mbox{ and }\\
&
|f(x,t)|\leq c_1(|t|^{\alpha} + 1),\ \forall t\in\mathbb{R}, \mbox{ for  a.e. } x \mbox{ in } \Omega,\\
& \mbox{for }1\leq \alpha<p^*-1 \mbox{ and } c_1>0,
\end{cases}
\end{equation}
and otherwise, for $p \geq N$, we suppose that
\begin{equation}
\begin{cases}
\tag{$H2^\prime$}\label{H2P}
&\mu\bigg(\{x\in \Omega: f(x,t)\neq t\}\bigg)>0 \mbox{ and }
\\
&|f(x,t)|\leq c_1(|t|^{\alpha} + 1),\ \forall t \in \mathbb{R} \mbox{ for a.e. } x \mbox{ in } \Omega,\\
& \mbox{for }1\leq \alpha<\infty \mbox{ and } c_1>0,
\end{cases}
\end{equation}
where $\mu$ stands for the $N$-dimensional Lebesgue measure.

Indeed, at this point, we must remark that, to the best of our knowledge, the assumption of 
\begin{equation}
\label{NT}\tag{$NT$}
\mu\bigg(\{x\in \Omega: f(x,t)\neq t\}\bigg)>0
\end{equation}
is new. Further, it is motivated by the fact that constant functions could possibly be solutions of \eqref{p1}, since we deal with Neumann boundary conditions. Thus, for the sake of convenience, we introduce the ``nontriviality condition" \eqref{NT}, which turns out to guarantee that constant solutions are avoided, usually an issue when dealing with this sort of boundary conditions.

Further, we stress that in order to obtain uniform $L^\infty$ a priori estimates for solutions of \eqref{p1}, we have to impose more restrictive conditions on $f(x,t)$. Actually, for this purpose, it is supposed that there exists $c_2>0$ and $\nu>0$ such that
\begin{equation}
\tag{$\mbox{NQ}$}
\label{NQ1}
\limsup_{|t|\to+\infty} \dfrac{f(x,t)t-\nu F(x,t)}{|t|^\sigma} \leq -c_2 \mbox{ uniformly for a.e. } x \mbox{ in } \Omega,
\end{equation}
where $F(x,t)=\displaystyle \int_0^t f(x,\tau) d\tau$ and $\sigma \in \mathbb{R}$.

For $\nu=p$, this is a version of the nonquadraticity hypothesis, originally introduced by Costa \& Magalh\~aes, for instance, see \cite{cm1}, which in the standard approach is employed to guarantee compactness conditions of Palais--Smale or Cerami type to a given functional whose critical points correspond to solutions of a given Partial Differential Equation. This is usually done by means of Variational Methods combined with an application of one of the several versions of the Mountain Pass theorem. Instead, in our case, we use \eqref{NQ1} 
 with a slightly distinct purpose -- it assures uniform a priori estimates for the solutions of \eqref{p1} -- which to the best of our knowledge is new. However, regardless of this usage being new, we stress that this choice is due to the sake of simplicity -- our purpose is to exploit fractional regularity methods for the investigation of nonlocal degenerate problems, see the remarks at the end of this section. 

Moreover, we are going to suppose that there exist $\beta \in [0,p^*-1)$ and $c_3>0$ such that
\begin{equation}
\tag{H3}\label{H3}
\limsup_{|t|\to+\infty} \dfrac{|F(x,t)|}{|t|^\beta} \leq c_3 <+\infty \mbox{ uniformly for a.e. } x \mbox{ in } \Omega,
\end{equation}
or
\begin{equation}
\tag{$H3^\prime$}\label{H3P}
\lim_{|t|\to+\infty} \dfrac{F(x,t)}{|t|^\beta} = -\infty \mbox{ uniformly for a.e. } x \mbox{ in } \Omega.
\end{equation}

Remark that condition \eqref{H3} means that $F(.,.)$ grows as  $|x|^\beta$ at infinity. Of course, this is a complementary and more accurate version of \eqref{H2}, since $\beta$ can always be assumed to be equal to $\alpha+1$. Nevertheless, the interesting case is to explore the possibility of having $\beta<\alpha+1$, what is absolutely not restrictive when compared to \eqref{H2}. Actually, if in the worst scenario \eqref{H3} holds for $\alpha=\beta+1$, since $\alpha <p^*-1$, one can always consider $\overline{\alpha}>\alpha$ such that $\overline{\alpha} <p^*-1$ and then switch $\alpha$ with $\overline{\alpha}$ in \eqref{H2} and then rename $\beta$ as $\alpha+1$. Further, condition \eqref{H3P} is stronger than \eqref{H2}, what allows us to obtain solutions for a broad class of exponents.

Now, we are in the position to introduce the main contributions of this work. In turn,  our main results address the existence and fractional regularity of a nontrivial strong solution for our problem, i.e., a nonconstant function satisfying \eqref{p1} a.e. in $\Omega$.

\begin{theorem}
\label{thm1}
Under hypotheses \eqref{H1}, \eqref{H2} or \eqref{H2P}, and \eqref{H3}, suppose that \eqref{NQ1} holds true. Then, for $p>2$, if
 \begin{equation}
 \nonumber
	 \nu c_3 < c_2 \mbox{ and } \sigma \geq \beta;
 \end{equation}
  there exists \[u \in \mathcal{N}^{1+\frac{2}{r_s},r_s}\setminus \{\mathbb{R}\},\] satisfying \eqref{p1} a.e. in $\Omega$,
where $\mathcal{N}^{1+\frac{2}{r_s},r_s}\setminus \{\mathbb{R}\}$ stands for the nonconstant elements of $\mathcal{N}^{1+\frac{2}{r_s},r_s}$.
 Moreover, there exists $C>0$ such that for all $s \in (2,+\infty)$,
\begin{equation}
\nonumber
\|u\|_{L^\infty}+\|u\|_{\mathcal{N}^{1+\frac{2}{r_s},r_s}}\leq K,
\end{equation}
where $C(N,p,s, a_0,c_1,c_2,c_3,\|a\|_{L^\infty_{loc}},\|F\|_{L^\infty_{loc}},\alpha,\beta,\sigma,\Omega)>0$.
\end{theorem}
By means of another choice for the assumptions on the primitive $F(.,.)$, we obtain our second result.
\begin{theorem}
\label{thm2}
Under hypotheses \eqref{H1}, \eqref{H2} or \eqref{H2P}, and \eqref{H3P}, suppose that \eqref{NQ1} holds true. Then, for $p>2$, we obtain the same conclusion of Theorem \ref{thm1}.
\end{theorem}

Concerning past results which are somehow connected to Theorem \ref{thm1} and Theorem \ref{thm2}, regardless of the specific boundary conditions and with no intention of being complete, we refer the reader to \cite{cog1,coru,rak,wxz} where $p$-Kirchhoff problems were investigated by the use of Variational Methods techniques. For instance, in \cite{cog1},  the authors investigate \eqref{p1} with Dirichlet boundary conditions where, among other assumptions, $f(.)$ and $a(.)$ are admitted to satisfy Ambrosetti--Rabinowitz conditions, and $f$ has subcritical growth. By using a version of the Mountain Pass theorem it is proved that \eqref{p1} possesses a positive solution.
In the case of Neumann boundary conditions, in \cite{coru}, the authors prove the existence of weak solutions to a system of $p$-Kirchhoff equations by means of the Ekeland Variational Principle, where it is assumed \eqref{H1} combined with certain smoothness and translation invariant hypotheses for $f$. More recently, in \cite{rak}, the threshold case $a(t)=a_1t+a_0$ was addressed with an approach based on a constraint variational method and an improved deformation lemma. By assuming certain growth and sign conditions for $f$, the authors prove that the Dirichlet version of \eqref{p1} admits at least one sing--changing solution. In addition, in \cite{wxz} existence and multiplicity results are proved as a consequence of a version of the Mountain Pass theorem, provided by the authors. In order to do that, it is assumed that $f$ and $a$ satisfy Ambrosetti--Rabinowitz conditions, and additional $p - q$ growth conditions for $f$, where $\inf\{a(t)\}$ is allowed to be zero.

Our technique to prove Theorems \ref{thm1}  and \ref{thm2} deeply relies on the derivation of explicit and global fractional order estimates for the solutions of \eqref{p1}, the aforementioned fractional regularity approach, and then, to employ the \eqref{NQ1} condition in order to obtain uniformly bounded solutions, see Proposition \ref{prop5}. This includes the investigation of a linearized version of \eqref{p1}, cf. Proposition \ref{prop1}, which generalizes the results obtained in \cite{dmm,dmm2} to nonlocal problems. As it turns out, the Kirchhoff term modifies the nature of the estimates in the cases where $p$ is close enough to $2$, see Remark \ref{remark1}.

We stress that for all of our estimates, a possible dependence on the norm of the solutions will be explicitly exhibited, mostly in terms of the $h_s$ polynomials, see Section \ref{apriori}. Actually, the a priori estimates given in Proposition \ref{prop1}, Lemma \ref{lemma2} and Propositions \ref{prop2}-\ref{prop4} provide explicit bounds for the norms of $\|u\|_{L^s}$ and $\|u\|_{\mathcal{N}^{1+\frac{2}{r_s},r_s}}$ depending on $h_s(\|u\|_{W^{1,p}})$ and $a(\|u\|^p_{W^{1,p}})$,  and have interest by their own since they may be adequate for other applications related to Kirchhoff problems. In turn, we remark that the aforesaid estimates, with the exception of the $L^\infty$ bounds given in Proposition \ref{prop5}, are independent of the \eqref{NQ1} condition. Indeed, \eqref{NQ1}  could be replaced by any of the known conditions which guarantee uniform $L^\infty$ bounds, since our concern in the present contribution relies on the use of the fractional regularity to relax the most common hypotheses on Kirchhoff--like terms. Actually, the only purpose of \eqref{NQ1}, \eqref{H3}, and \eqref{H3P} is to guarantee the validity of uniform energy bounds for the solutions.

Now, for the sake of completeness, we provide two examples where the conditions of the latter theorems hold.

\textbf{Example 1.}
First, we address the case where the conditions of Theorem \ref{thm1} hold. Thence, let us fix $\beta=\sigma \in(1,p]$ and set $f(x,t)=-ct(t^2+1)^{\frac{\beta-2}{2}}+g(x)$, where $g\in L^\infty \setminus \{\mathbb{R}\}$ and $c>0$. Moreover, consider
\[a(t)=\theta_1\ln(1+|t|)+\theta_2, \]
for $\theta_i>0$, i=1, 2.
Once again, observe that $a(.)$ cannot be bounded by below by any sort of polynomial growth  and does not satisfy the Ambrosetti--Rabinowitz conditions.

In addition, it is obvious that $a(.)$ and $f(.,.)$ respectively satisfy \eqref{H1} and \eqref{H2}. Complementarily, in this case the primitive $F(.,.)$ is given by
\[F(x,t)=-\dfrac{c(t^2+1)^{\frac{\beta}{2}}+c}{\beta}-tg(x),
\]
so that \eqref{H3} and \eqref{NQ1} are satisfied if $\nu <\frac{\beta}{2}$. We stress that $c_2=\frac{c(\beta-\nu)}{\beta}$ and $c_3=\frac{c 2^{\beta/2}}{\beta}$ so that by the choice of $\nu$, we get $\nu c_3<c_2$, and thus Theorem \ref{thm1} is applicable.

\textbf{Example 2.}
Now, we present an example where the conditions of Theorem \ref{thm2}, hold. In turn, let
\begin{align*}
f(x,t)&=-g(x)t\bigg[(p-2)(p-\epsilon)|t|^{p-\epsilon-2}\sin^2\left(\frac{|t|^{\epsilon}}{\epsilon}\right) \\
& + \left(p+(p-2)\sin\left(\frac{2|t|^{\epsilon}}{\epsilon}\right)\right)|t|^{p-2}\bigg],
\end{align*}
c.f. \cite{ChMa}, where $p>2$, $0<\epsilon <p-2$, $g\in L^\infty \setminus \mathbb{R}$, $g\geq g_0>0$ with $g_0$ a constant,
so that $F(x,t)=-g(x)\left(|t|^p+(p-2)|t|^{p-\epsilon}\sin^2\left(\frac{|t|^{\epsilon}}{\epsilon}\right)\right)$. Then consider
\[a(t)=
\begin{cases}
 &\delta_1\bigg|t \sin{\bigg(\dfrac{1}{t}\bigg)}\bigg|+\delta_2 \mbox{ if } 0<t \leq 1,\\
 &\delta_1+\delta_2 \mbox{ if } t= 0,
 \end{cases}
 \]
for $\delta_i>0$, i=1,2, and where $a(.)$ is defined in the rest of $\mathbb{R}$ by means of a periodic extension.
Remark that due to its highly oscillatory pattern, $a(.)$ is nonsmooth, actually, $a(.)$ is not even in $BV$. Further, it cannot be bounded by below by any sort of polynomial growth.

Now, it is clear that $a(.)$ and $f(.,.)$ respectively satisfy \eqref{H1} and \eqref{H2}. Moreover, $f(.,.)$ and $F(.,.)$ satisfy \eqref{NQ1} and \eqref{H3P}, for \[0<\beta<\min\{p-\epsilon, p^*-1\} \mbox{ and } \sigma=p.\] Indeed, let $0<\theta,\nu$ such that
\[
\frac{1}{2}+\theta=\frac{1}{\nu},
\]
that is, $\nu=\frac{2}{1+2\theta}$, then it is clear that
\[
\begin{array}{l}
\frac{1}{\nu}tf(x,t)\\
=\frac{-g(x)}{2}\left[(p-2)(p-\epsilon)|t|^{p-\epsilon}\sin^2\left(\frac{|t|^{\epsilon}}{\epsilon}\right) + \left(p+(p-2)\sin\left(\frac{2|t|^{\epsilon}}{\epsilon}\right)\right)|t|^{p}\right]\\
-\theta g(x)\left[(p-2)(p-\epsilon)|t|^{p-\epsilon}\sin^2\left(\frac{|t|^{\epsilon}}{\epsilon}\right) + \left(p+(p-2)\sin\left(\frac{2|t|^{\epsilon}}{\epsilon}\right)\right)|t|^{p}\right]\\
\leq -g(x)\left[(p-2)|t|^{p-\epsilon}\sin^2\left(\frac{|t|^{\epsilon}}{\epsilon}\right) + |t|^{p}\right]-\theta g_0|t|^{p}\\
= F(x,t) - \nu\theta g_0|t|^{p}.
\end{array}
\]
Hence,
\[
\frac{tf(x,t) -\nu F(x,t)}{|t|^{p}}\leq - \nu\theta g_0 <0, \forall t \in \mathbb{R}\setminus \{0\},
\]
and \eqref{NQ1} holds.

Further, once
\[
\frac{F(x,t)}{|t|^{\beta}}=-g(x)\left(|t|^{p-\beta}+(p-2)|t|^{p-\epsilon-\beta}\sin^2\left(\frac{|t|^{\epsilon}}{\epsilon}\right)\right),
\]
thus \eqref{H3P} also holds if $0<\beta<\min\{p-\epsilon, p^*-1\}$. And then, Theorem \ref{thm2} is applicable in this case.

Finally, we stress that throughout the paper, the symbol $C$ denotes a general constant which does not depend on the solutions of \eqref{p1}, $u$, and may vary from line to line. Moreover, sometimes we are going to denote that $C$ depends on a function which is already defined in terms of the datum. This means that $C$ blows up at the singularities of the given function. For instance, $C=C(N,p,p^*)$ means that this constant degenerates whenever $p^*=+\infty$, and so on. Further, the dependence on $\Omega$ of the functions spaces considered in the text, unless otherwise stated, is going to be omitted.


\section{Preliminary results}
\label{PreliminaryResults}
In this section, we investigate a linearized version of \eqref{p1} and obtain the preliminary fractional energy estimates employed along the text.
We start with an algebraic iterative lemma which, for the convenience of the reader, will have its proof exhibited.
\begin{lemma}\label{algel} Consider $\{A_n\}, \{d_n\}$, sequences of real numbers, and $\kappa\in \mathbb{R}$. Suppose that
\begin{equation}
\nonumber
A_n\leq C (1+A_{n-1}+A_{n-1}^{d_{n-1}}), A_n \geq 0, d_n \geq 0, d_n \leq d_{n+1}  \forall n \in \mathbb{N},
\end{equation}
and that $d_n \to \kappa$, if $n\to \infty$. Then, there exists $K_n$ such that
\begin{equation}
\label{2algel}
A_n \leq K_n\big(1+A_0+A_0^{n \cdot \kappa}\big), \forall n \in \mathbb{N},
\end{equation}
where $K_n=K_n(C,\kappa)>0$ and $K(n)\to +\infty$, if $n\to +\infty$.
\end{lemma}
\begin{proof}

Suppose that inequality \eqref{2algel} holds for all $n \leq m-1$. By hypothesis, we have
\begin{align}
\nonumber
A_m &\leq C (1+A_{m-1}+A_{m-1}^{d_{m-1}})
\\
\nonumber
&\leq C\bigg( 1+K_{m-1}\big(1+A_0+\big(A_0^d\big)^{m-1}\big)
+K_{m-1}^{d_{m-1}}\big(1+A_0+\big(A_0^d\big)^{m-1}\big)^{d_{m-1}}\bigg)
\\
\nonumber
&\leq K_m\bigg(1+A_0+\big(A_0^d\big)^m\bigg).
\end{align}
\end{proof}



Our next proposition regards the existence of solutions and a priori bounds in Nikolskii spaces, $\mathcal{N}^{s,t}$, for an auxiliary version of \eqref{p1}.
This result, which has interest by its own, improves to the context of nonlocal problems the fractional regularity obtained in \cite{dmm,dmm2}.
\begin{proposition}
\label{prop1} Suppose that \eqref{H1} holds. Then, given $g\in L^s$, for $s>2$, there exists $u \in \mathcal{N}^{1+\frac{2}{r_s},r_s}$ such that
\begin{equation}
\nonumber
\left\{
\begin{aligned}
-a(\|u\|^p_{W^{1,p}})^{p-1}\Delta_p u + u & =  g && \text{a.e. in } \Omega,\\
     \dfrac{\partial u}{\partial\eta} & =    0 && \text{on } \partial \Omega,
\end{aligned}
\right.
\end{equation}
where  $r_s= s(p-2)+2$, see Identity \eqref{exponents2}.

Moreover,

\begin{equation}
\label{pr1}
\|u\|_{\mathcal{N}^{1+\frac{2}{r_s},r_s}} \leq C\bigg(1+\|g\|_{L^s}+a\big(\|u\|^p_{W^{1,p}}\big)^{w_{p,s}} \|g\|^{s/r_s}_{L^s}\bigg),
\end{equation}
where $C=C(N, p, a_0, s, \Omega)>0$.
\end{proposition}

\begin{proof}

The idea is to combine the fractional regularity results given in \cite{dmm,dmm2}, with the Leray--Schauder Fixed Point Theorem and a scaling technique.

Indeed, for the purpose of this proof, let us consider the parameter $t\in[0,1]$ and then set the following family of operators:
\[T_t:W^{1,p}\to W^{1,p} \mbox{ for } 0\leq t \leq 1, \mbox{ where } T_t(v)=u \mbox{ if and only if }\]

\begin{equation}
\label{pr2}
\left\{
\begin{aligned}
-b(v)\Delta_p u + u & =  t g && \text{in } \Omega,\\
     \dfrac{\partial u}{\partial\eta} & =    0 && \text{on } \partial \Omega.
\end{aligned}
\right.
\end{equation}

We will check the hypotheses of the Leray--Schauder in five steps.

\noindent \textit{Step 1.} $T_t$ is well--defined for every $t\in[0,1]$ and $T_t(\wpp) \hookrightarrow\hookrightarrow \wpp$.

Indeed, given $(v,t)\in \wpp \times [0,1]$, set $k=(b(v))^{1/(p-2)}$ and $\hat{g}=t kg \in L^s$. Further, by \cite{dmm}, see Thm. 1.1, and \cite{dmm2}, see Thm. 2.2, there exists $\hat{u} \in \mathcal{N}^{1+\frac{2}{r_s},r_s}\cap \wpp$, the {\bf unique} strong solution of
\begin{equation}
\label{pr3}
\left\{
\begin{aligned}
-\Delta_p \hat{u} + \hat{u} & =  \hat{g} && \text{a.e. in } \Omega,\\
     \dfrac{\partial \hat{u}}{\partial\eta} & =    0 && \text{on } \partial \Omega.
\end{aligned}
\right.
\end{equation}
 which satisfies the following fractional energy estimate

\begin{equation}
\nonumber
\|\hat{u}\|_{\mathcal{N}^{1+\frac{2}{r_s},r_s}} \leq C\bigg(1+\|\hat{g}\|_{L^s}+\|\hat{g}\|_{L^s,}^{s/r_s}\bigg),
\end{equation}
where $C=C(N, p, s, \Omega)>0$.

Now, it is convenient to consider $u=\hat{u}/k \in \mathcal{N}^{1+\frac{2}{r_s},r_s}.$ Then, since $k^{p-2}=b(v)$, there holds that

\begin{equation}
\label{pr5}
\left\{
\begin{aligned}
-b(v)\Delta_p u + u & =  t g && \text{a.e. in } \Omega,\\
     \dfrac{\partial u}{\partial\eta} & =    0 && \text{on } \partial \Omega.
\end{aligned}
\right.
\end{equation}
As well as there exists a unique solution to \eqref{pr3}, the solution of \eqref{pr5} is also unique, and thenceforth $T_t$ is well--defined.

In order to prove $T_t$'s compactness, first observe that, since $p>2$,
\[
k^{2/p-1}=k^{(2-p)/p} \leq a_0^{(2-p)/p} \mbox{ and } k^{s/r_s-1} \leq (1+a_0^{s/r_s-1})a(\|v\|^p_{W^{1,p}})^{\omega_{p,s}},
\]
where the exponent $w_{p,s}$ is given in \eqref{exponents4}, p. \pageref{exponents4}.

Then, by combining \eqref{pr3}, the definition of $\hat{u}$ and by performing some simplifications, we end up with
\begin{align}
\nonumber
\|u\|_{\mathcal{N}^{1+\frac{2}{r_s},r_s}} &\leq C\bigg(\frac{1}{k}+t \|g\|_{L^s}+k^{s/r_s-1}t^{s/r_s}\|g\|_{L^s}^{s/r_s}\bigg)
\\
\label{pr6}
&\leq C\bigg(a_0^{\frac{1-p}{p-2}}+\|g\|_{L^s}+a(\|v\|^p_{W^{1,p}})^{\omega_{p,s}}\|g\|_{L^s}^{s/r_s}\bigg).
\end{align}
In addition, since $r_s>p$, remark that $\mathcal{N}^{1+\frac{2}{r_s},r_s} \hookrightarrow\hookrightarrow \wpp$ and then, by \eqref{pr6} we see that \[T_t(\wpp)\hookrightarrow\hookrightarrow \wpp, \mbox{ for all } t \in [0,1].\]

\noindent \textit{Step 2.} For every $t\in [0,1]$  fixed,  $T_t:\wpp \to \wpp$ is continuous.

Consider $\{v_n\}_{\mathbb{N}} \subset \wpp$ such that $v_n \to v$ in $\wpp$, if $n \to \infty$. Moreover,  set $u_n = T_t(v_n)$ and $u=T_t(v)$.

By \eqref{pr6}, it follows that

\begin{align*}
\|u_n\|_{\mathcal{N}^{1+\frac{2}{r_s},r_s}}
&\leq C\bigg(1+\|g\|_{L^s}+a(\|v_n\|^p_{W^{1,p}})^{\omega_{p,s}}\|g\|_{L^s}^{s/r}\bigg),
\\
&\leq C\bigg(1+\|g\|_{L^s}+\|g\|_{L^s}^{s/r}\bigg) \forall n \in \mathbb{N},
\end{align*}

since $a(.)$ is continuous.

Then, there exists $\hat{u} \in \mathcal{N}^{1+\frac{2}{r_s},s}$ such that, up to subsequences,
\begin{equation}
\label{pr7}
u_n\rightharpoonup \hat{u} \mbox{ in } \mathcal{N}^{1+\frac{2}{r_s},s} \hookrightarrow\hookrightarrow \wpp.
\end{equation}
Now, by a standard argument, it is enough to prove that $\hat{u}=u$, i.e., $\hat{u}$ is the weak solution of

\begin{equation}
\label{pr8}
\left\{
\begin{aligned}
-b(v)\Delta_p \hat{u} + \hat{u} & =  t g && \text{in } \Omega,\\
     \dfrac{\partial \hat{u}}{\partial\eta} & =    0 && \text{on } \partial \Omega.
\end{aligned}
\right.
\end{equation}

For this, remark that for every $\phi \in \wpp$ we have

\begin{align}
\nonumber
& \ \ \bigg|\int_\Omega b(u_n)|\nabla u_n|^{p-2}\nabla u_n \cdot \nabla \phi - \int_\Omega b(v) |\nabla \hat{u}|^{p-2}\nabla \hat{u} \cdot \nabla \phi\bigg|
 \\
\nonumber
&\leq |b(v_n)-b(v)| \|u_n\|^{p-1}_\wpp \|\phi\|_\wpp
+b(v)\int_\Omega\bigg||\nabla u_n|^{p-2}\nabla u_n - |\nabla \hat{u}|^{p-2}\nabla \hat{u}| |\nabla \phi|
\\
\label{pr9}
&C|b(v_n)-b(v)|\|\phi\|_\wpp+b(v)\int_\Omega\big||\nabla u_n|^{p-2}\nabla u_n - |\nabla \hat{u}|^{p-2}\nabla \hat{u}\big| |\nabla \phi|.
 \end{align}

In addition, remark that
\begin{equation}
\nonumber
  |b(v_n)-b(v)|=\bigg| a\big(\|v_n\|^p_\wpp)^{p-1}-a\big(\|v\|^p_\wpp)^{p-1}\bigg|\to 0 \mbox{ as } n \to \infty
\end{equation}
since $a(.)$ is continuous and $v_n \to v$ in $\wpp$.

Moreover, by \eqref{pr7} there exists another subsequence, relabelled the same, such that
\[ \nabla u_n \to \nabla \hat{u} \mbox{ a.e. in } \Omega \mbox{ and } |\nabla u_n| \leq |\Gamma|,
\]
for a fixed function $\Gamma \in L^p$.

Consequently,

\[ |\nabla u_n|^{p-2}\nabla u_n-|\nabla \hat{u}|^{p-2}\nabla \hat{u} \to 0 \mbox{ a.e. in }\Omega,\]
and

\[ \big||\nabla u_n|^{p-2}\nabla u_n-|\nabla \hat{u}|^{p-2}\nabla \hat{u} \big| |\nabla \phi| \leq \big(|\omega|^{p-1}+|\nabla \hat{ u}|^{p-1}\big)|\nabla \phi| \in L^1.
\]

Hence, by the Lebesgue Dominated Convergence theorem

\begin{equation}
\label{pr11}
\int_\Omega\big||\nabla u_n|^{p-2}\nabla u_n - |\nabla \hat{u}|^{p-2}\nabla \hat{u}\big| |\nabla \phi|\to 0.
\end{equation}

Thus, by combining \eqref{pr9}-\eqref{pr11} we deduce that $\hat{u}$ is a weak solution of \eqref{pr8}, which has a unique solution by Step 1. Thence,  $\hat{u}=u$ a.e. in $\Omega$, and by a standard argument, we see that $T_t$ is continuous.

\noindent \textit{Step 3.} $T_t(.)$ is uniformly continuous with respect to $t$.

Now, take $u_i=T_{t_i}(v)$, for $i=1, 2$  and set $u=u_1-u_2$. Since
\[-b(v)\Delta_p u_i +u_i=t_i g,\]
by taking the difference between the cases $i=1$ and $i=2$, by multiplying the result by $u$, by integrating by parts, owing to Tartar's (see  Lemma 4.4, \cite{dibenedetto} p.14 ) and Young's inequalities, there follows that

\[\int_\Omega|\nabla u|^p+\int_\Omega |u|^2 \leq C|t_1-t_2|\|g\|^2_{L^s}.\]

In this way, by Gagliardo--Nirenberg's interpolation inequality, we get
\[\|u\|_\wpp \leq C|t_1-t_2|^{1/p}\|g\|_{L^s}^{2/p} + C|t_1-t_2|^{1/2}\|g\|_{L^s},\]
and this proves that $T_t$ is indeed uniformly continuous with respect to $t$.

\noindent \textit{Step 4.} There exists $C>0$ such that for every fixed point of $T_1$
\[\|u\|_\wpp \leq C.\]
Now,  by recalling \eqref{pr6}, given $u=T_1(u)$ we arrive at
\begin{equation}
\label{pr12}
\|u\|_{\mathcal{N}^{1+\frac{2}{r_s},r_s}} \leq C \bigg(1+\|g\|_{L^s}+a(\|u\|^p_{W^{1,p}})^{\omega_{p,s}}\|g\|_{L^s}^{s/r_s}\bigg),
\end{equation}
which is bounded by $C (1+\|g\|_{L^s}+\|g\|_{L^s}^{s/r_s})$, where $C=C(N,p,s,a_0,\Omega)>0$, since $a(.)$ and continuous. The latter estimate finishes the proof of Step 4.

\noindent \textit{Step 5.} Finally, it remains to check that $T_0(.,v)\equiv 0$. This is a direct consequence of Definition \eqref{pr2}.

Hence, by Steps 1-5 and by the Leray--Schauder Fixed Point Theorem, there exists $u\in \wpp$ such that $T_1 u=u$, or equivalently

\begin{equation}
\nonumber
\left\{
\begin{aligned}
-b(u)\Delta_p u + u & =  g && \text{in } \Omega,\\
     \dfrac{\partial u}{\partial\eta} & =    0 && \text{on } \partial \Omega.
\end{aligned}
\right.
\end{equation}

Moreover, by the very definition of $T_t$,  $u\in \mathcal{N}^{1+\frac{2}{r_s},r_s}$ satisfies
\[-a(\|u\|^p_{W^{1,p}})^{p-1}\Delta_p u + u  =  g \text{ a.e. in } \Omega,\]
see the proof of Theorem 2.2 in \cite{dmm2}. Actually, since
\[|\nabla u|^{p-2}\nabla u \in L^{p^\prime}\]
and, by now, $b(u)$ could merely be considered as a constant, by applying Theorem III.2.2 from \cite{galdi}, $u$ fulfills
\[\int_\Omega - b(u)\Delta_p u  \phi +u \phi = \int_\Omega f \phi, \ \  \forall \phi \in W^{1,p},\]
and thus it is indeed a strong solution.

Finally, remark that from inequality \eqref{pr12}, $u$ also satisfies \eqref{pr1}, what concludes the proof of Proposition \ref{prop1}.

\end{proof}

\begin{remark}
\label{remark1}
We stress that based on estimate \eqref{pr1}, the Kirchhoff term favors fractional order energy estimates. 
 Moreover, observe that when $p>3-\frac{2}{s}$ the effects of this term disappear.
\end{remark}
\section{A priori estimates}
\label{apriori}
In this section, we establish  explicit a priori bounds for solutions of \eqref{p1} with respect to $L^\infty$ and fractional order spaces.  First, we obtain estimates which depend on appropriate norms of the solutions of \eqref{p1}. Later on, by using condition \eqref{NQ1},
 we provide uniform a priori bounds for the solutions.
 We stress that in the following results, we are going to employ the definition of $h_s$ polynomials, see \eqref{2lm2} on p. \pageref{2lm2} and the related notations. Nevertheless, these estimates are new and provide a precise measure of the interplay between the degeneracy parameter $p$, the Kirchhoff term $a(.)$, the nonlinearity $f(.,.)$, and the space dimension $N$.

As a first step, we investigate the case $2<p<N$ for $\alpha \geq p-1$, and by following the classical argument due to J. Moser, we prove that every solution of \eqref{p1} bootstraps itself into improved Lebesgue spaces.

\begin{lemma} \label{lemma2} Under hypotheses \eqref{H1} and \eqref{H2}, suppose that $2<p<N$ and $\alpha \geq p-1$. Given $u\in W^{1,p}$, a weak solution of \eqref{p1}, consider $\delta>0$ such that $\alpha=p^*-1-\delta p^*$. Then, given $s\in (2,+\infty)$, there holds that

\begin{equation}
\label{1lm2}
\|u\|_{L^s} \leq C \bigg(1+h_s\big(\|u\|_{L^{p^*}}\big)\bigg)
\end{equation}

where $C=C(N, p, p^*, n_s, a_0,c_1, \alpha, \delta, \Omega)>0$ and 
\[n_s=\min\{n\in \mathbb{N}:p^* (k_n+1)\geq s\}, \mbox{ see \eqref{4lm2}, }\]
for $k_n= \delta \sum_{i=1}^{n+1} \lambda^i$.

\end{lemma}
\proof

Set $k\geq0$. The proof consists in using $v_m^{kp+1}$ as a test function in \eqref{p1}, where
\[v_m=\min\{u_+,m\} \in W^{1,p}\cap L^\infty, u_+=\max\{u,0\} \mbox{ and }m>0,\]
in order to obtain certain iterative estimates for $u$.

For the sake of clarity, we split the argument into four steps.

\noindent \textit{Step 1.} {\it Higher order estimates}.

Owing to standard calculations, it is true that
\begin{align*}
&\int_\Omega b(u)|\nabla u|^{p-2}\nabla u \cdot \nabla v_m^{kp+1}+u v_m^{kp+1}&
\\
&> \dfrac{a_0^{p-1}(kp+1)}{(k+1)^p}\int_\Omega |\nabla v_m^{k+1}|^p+ \int_\Omega v_m^{kp+2}.&
\end{align*}
Then, there holds the following estimate
\begin{align}
\nonumber
&\int_\Omega b(u)|\nabla u|^{p-2}\nabla u \cdot \nabla v_m^{kp+1}+u v_m^{kp+1}
\\
\nonumber
&>\dfrac{a_0^{p-1}(kp+1)}{(k+1)^p}\|\nabla v_m^{k+1}\|^p_{L^p}+ \|v_m\|^{kp+2}_{L^{kp+2}}, \forall m>0 \mbox{ and } k\geq0.
\end{align}

\noindent \textit{Step 2.} {\it Lower order estimates}.

By hypothesis \eqref{H2},
\[\int_\Omega f(x,u) v_m^{kp+1} \leq c_1 \int_\Omega (|u|^\alpha+1)v_m^{kp+1}.\]

Then, remark that
\begin{align}
\nonumber
\int_\Omega |u|^\alpha v_m^{kp+1} &= \int_\Omega u_+^\alpha v_m^{kp+1}
\\\nonumber
&\leq \int_\Omega u_+^{\alpha+kp+1}=\|u_+\|^{\alpha+kp+1}_{L^{\alpha+kp+1}}.
\end{align}

Moreover, by H\"{o}lder's inequality

\[
\int_\Omega v_m^{kp+1} \leq |\Omega|^{\alpha/(\alpha+kp+1)}\|v_m\|^{kp+1}_{L^{\alpha +kp+1}}\leq C\|u_+\|^{kp+1}_{L^{\alpha +kp+1}},
\]
since $|\Omega|^{\alpha/(\alpha+kp+1)}\leq |\Omega|+1$.

In this fashion, the following estimate holds true
\begin{align*}
\int_\Omega f(x,u)v_m^{kp+1}&\leq C\bigg(\|u_+\|^{kp+1}_{L^{\alpha + kp+1}}+\|u_+\|^{\alpha+kp+1}_{L^{\alpha+kp+1}}\bigg)
\end{align*}
where $C=C(c_1, \alpha, \Omega)>0$.

\noindent \textit{Step 3.} {\it Iterative estimates}

Recalling that
\[
\int_\Omega a(u) |\nabla u|^{p-2}\nabla u \cdot \nabla \phi +u \phi = \int_\Omega f(x,u)\phi, \forall \phi \in W^{1,p},
\]
we fix $\phi=v_m^{kp+1}$ and combine Steps 1 and 2, obtaining in this manner
\begin{align*}
&\dfrac{a_0^{p-1}(kp+1)}{(k+1)^p}\|\nabla v_m^{k+1}\|^p_{L^p}+ \|v_m\|^{kp+2}_{L^{kp+2}}
\\
&< C\bigg(\|u_+\|^{kp+1}_{L^{\alpha+kp+1}}+\|u_+\|^{\alpha+kp+1}_{L^{\alpha+kp+1}}\bigg), \forall m>0 \mbox{ and } k\geq0.
\end{align*}

Then, owing to the Fatou Lemma, by letting $m\to+\infty$, there follows
\begin{align*}
&\dfrac{a_0^{p-1}(kp+1)}{(k+1)^p}\|\nabla u_+^{k+1}\|^p_{L^p}+ \|u_+\|^{kp+2}_{L^{kp+2}}
\\
&< C\bigg(\|u_+\|^{kp+1}_{L^{\alpha+kp+1}}+\|u_+\|^{\alpha+kp+1}_{L^{\alpha+kp+1}}\bigg), \forall k\geq0,
\end{align*}
for $v_m\to u_+$ and $\nabla v_m \to \nabla u_+$ a.e. in $\Omega$.

However, remark that adding on both sides of the last inequality the term $\displaystyle\|u_+^{k+1}\|^p_{L^p}$, leads to
\begin{align*}
&\dfrac{a_0^{p-1}(kp+1)}{(k+1)^p}\|\nabla u_+^{k+1}\|^p_{L^p}+ \|u_+^{k+1}\|^{p}_{L^{p}}
\\
&< C\bigg(\|u_+\|^{kp+1}_{L^{\alpha+ kp+1}}+\|u_+\|^{\alpha+kp+1}_{L^{\alpha+kp+1}}+\|u_+\|^{(k+1)p}_{L^{(k+1)p}}\bigg).
\end{align*}

Thus, by combining Bernoulli's inequality
\[
kp+1\leq (k+1)^p, \forall k\geq0,
\]
 and Sobolev's Embedding theorem, we arrive at
\begin{align}
\nonumber
&\dfrac{(kp+1)}{(k+1)^p}\|u_+^{k+1}\|^p_{L^{p^*}}
\\
\nonumber
&< C\bigg(\|u_+\|^{kp+1}_{L^{\alpha+kp+1}}+\|u_+\|^{\alpha+kp+1}_{L^{\alpha+kp+1}}+\|u_+\|^{(k+1)p}_{L^{(k+1)p}}\bigg),
\end{align}
where $C=C(N, p, a_0, c_1, \alpha, \Omega)>0$.
Nevertheless, for
\[
 \|u_+\|_{L^{(k+1)p^*}}=\|u_+^{k+1}\|_{L^{p^*}}^{1/(k+1)},
 \]
there holds that
\begin{align}
\nonumber
\|u_+\|_{L^{p^*(k+1)}}&< \bigg[C\dfrac{(k+1)^p}{(kp+1)}\bigg(\|u_+\|^{kp+1}_{L^{\alpha+kp+1}}
\\
\label{5lm2}
&+\|u_+\|^{\alpha+kp+1}_{L^{\alpha+kp+1}}+\|u_+\|^{(k+1)p}_{L^{(k+1)p}}\bigg)\bigg]^{1/p(k+1)}.
\end{align}
Now, we have to adjust the latter exponents. First, recall that we are considering $\alpha\geq p-1$, so that
\[
\int_\Omega u_+^{(k+1)p} \leq |\Omega|^{(\alpha-p+1)/(\alpha+kp+1)} \|u_+\|_{L^{\alpha+kp+1}}^{(k+1)p} \leq (|\Omega|+1)\|u_+\|_{L^{\alpha+kp+1}}^{(k+1)p}.
\]

Thus, by \eqref{5lm2} we are led to

\begin{align}
\nonumber
\|u_+\|_{L^{p^*(k+1)}}&<\bigg[C\dfrac{(k+1)^p}{(kp+1)}\bigg(\|u_+\|^{kp+1}_{L^{\alpha+kp+1}}
\\
\label{6lm2}
 &+\|u_+\|^{\alpha+kp+1}_{L^{\alpha+kp+1}}+\|u_+\|^{(k+1)p}_{L^{\alpha+kp+1}}\bigg)\bigg]^{1/p(k+1)}.
\end{align}

 In addition, there holds
\begin{equation}
\nonumber
(B_1+B_2+B_3)^{1/p(k+1)}\leq 2^{2/p(k+1)}(B_1^{1/p(k+1)}+B_2^{1/p(k+1)}+B_3^{1/p(k+1)}),
\end{equation}
for all $B_i>0$, and also remark that
\begin{equation}
\nonumber
 C^{1/(k+1)}(k+1)^{1/(k+1)}\leq C e^{1/e}, \forall k\geq0.
\end{equation}
In this way, inequality \eqref{6lm2} leads to
\begin{align}
\nonumber
\|u_+\|_{L^{p^*(k+1)}}&<\bigg(\dfrac{2C(k+1)^p}{(kp+1)}\bigg)^{1/p(k+1)}\bigg(\|u_+\|^{(kp+1)/p(k+1)}_{L^{\alpha+kp+1}}
\\
\nonumber
&+\|u_+\|^{(\alpha+kp+1)/p(k+1)}_{L^{\alpha+kp+1}}+\|u_+\|_{L^{\alpha+kp+1}}\bigg)
\\
\nonumber
&
\leq C^{1/(k+1)}(k+1)^{1/(k+1)}\bigg(1+\|u_+\|_{L^{\alpha+kp+1}}
\\
\nonumber
&+\|u_+\|^{(\alpha+kp+1)/p(k+1)}_{L^{\alpha+kp+1}}\bigg)
\\
&
\leq C\bigg(1+\|u_+\|_{L^{\alpha+kp+1}}+\|u_+\|^{(\alpha+kp+1)/p(k+1)}_{L^{\alpha+kp+1}}\bigg),
\label{7lm2}
\end{align}
where $C=C(N, p, a_0, c_1, \alpha, \Omega)>0$.

Observe that, in particular, the last inequality guarantees that $u\in L^s$ for all $s>1$.  Indeed, recall that for all $k>0$, we have \[p^*(k+1)> \alpha +kp+1, \mbox{ since }p^*>\alpha+1,\] by \eqref{H2}.

Then, if we assume that there exists $\bar{s}=\sup\{s\geq 1: u \in L^s\}\in \mathbb{R}$, by setting $\bar{k}_\epsilon$ such that
\[\alpha+\bar{k}_\epsilon p+1=\bar{s}-\epsilon,\]
we arrive at a contradiction if we choose $\epsilon>0$, for which
\[p^*(\bar{k}_\epsilon+1)>\bar{s}.\]

\noindent \textit{Step 4.} {\it $L^s$ estimates}.

Consider $\delta>0$ and $k_0\geq0$, for which
\begin{align}
\nonumber \alpha &=(1-\delta)p^*-1
\\
\nonumber
k_0&=\lambda-\dfrac{\alpha+1}{p}=\delta \lambda \mbox{ , and }
\\
\label{sequence1}
k_n&= \delta \sum_{i=1}^{n+1} \lambda^i,
\end{align}
where $\lambda=\frac{p^*}{p}$, see \eqref{exponents3}.

By these choices, there holds that
\[
\alpha+k_0p+1=p^* \mbox{ and } \alpha+k_np+1=p^*(1+k_{n-1}), \forall n \in \mathbb{N}.
\]
Then, owing to inequality \eqref{7lm2}
\begin{align*}
\nonumber
&\|u_+\|_{L^{p^*(k_n+1)}}
\\
&
\leq C\bigg(1+\|u_+\|_{L^{\alpha+k_np+1}}+\|u_+\|^{(\alpha+k_np+1)/p(k_n+1)}_{L^{\alpha+k_np+1}}\bigg)
\\
&
= C\bigg(1+\|u_+\|_{L^{p^*(1+k_{n-1})}}+\|u_+\|^{\lambda(1+k_{n-1})/(1+k_n)}_{L^{p^*(1+k_{n-1})}}\bigg),
\end{align*}
where $C=C(N, p, p^*, a_0, c_1, \alpha, \Omega)>0$. In particular, for $n=0$, there holds that
\begin{align*}
\nonumber
\|u_+\|_{L^{p^*(k_0+1)}}
&\leq C\bigg(1+\|u_+\|_{L^{\alpha+k_0 p+1}}+\|u_+\|^{(\alpha+k_0 p+1)/p(k_0+1)}_{L^{\alpha+k_0p+1}}\bigg)
\\
&=C\bigg(1+\|u_+\|_{L^{p^*}}+\|u_+\|^{\lambda/(\delta \lambda +1)}_{L^{p^*}}\bigg),
\end{align*}
by the choice of the parameters $k_0$ and $\delta$.

Then, by combining the latter inequalities with Lemma \ref{algel}, for  $A_n=\|u_+\|_{L^{p^*(k_n+1)}}$, $d_n=\lambda(1+k_{n-1})/(1+k_n)$, and $\kappa=\lambda$, we arrive at
\begin{align*}
\|u_+\|_{L^{p^*(k_n+1)}} &\leq C\bigg (1+\|u_+\|_{L^{p^*(k_0+1)}}+ \|u_+\|_{L^{p^*(k_0+1)}}^{n \cdot \lambda}\bigg)\\
&\leq C\bigg (1+\|u_+\|_{L^{p^*}}+ \|u_+\|_{L^{p^*}}^{(n+1) \lambda}\bigg), \forall n \in \mathbb{N},
\end{align*}
where $C=C(N, p, p^*, n, a_0, c_1, \alpha, \delta, \Omega)>0$.

Now, given $s \in (1,+\infty)$, take
\[ n_s =\min_{n\in \mathbb{N}} \{p^*(k_n+1)\geq s\}.
\]
Thus,
\[
\|u_+\|_{L^s} \leq C\bigg (1+\|u_+\|_{L^{p^*}}+ \|u_+\|_{L^{p^*}}^{(n_s+1)\lambda}\bigg)<+\infty,
\]
where $C=C(N, p, p^*, n_s, a_0, c_1, \alpha, \delta, \Omega)>0$.

In an analogous manner, see the proof of Lemma \ref{lemma4} below, by considering
\[
v_m=\min\{u_-,m\}, \mbox{ where } u_-=\max\{-u,0\},
\]
we obtain the latter estimate for $u_-$, and therefore the result follows.


With minor modifications on the arguments, we prove that the same result holds for the simpler case $1\leq \alpha <p-1$.

\begin{lemma} \label{lemma4} Under hypotheses \eqref{H1} and \eqref{H2}, suppose that $2<p<N$ and $1\leq \alpha < p-1$. Given $u\in W^{1,p}$, a weak solution of \eqref{p1} and $s\in (2,+\infty)$, there holds that

\begin{equation}
\nonumber
\|u\|_{L^{s}}\leq C\bigg(1+\|u\|_{L^{p^*}}\bigg)
\end{equation}

where $C=C(N, p, p^*, n_s, a_0,c_1, \alpha, \Omega)>0$ and $n_s=\min\{n\in \mathbb{N}:p^* (k_n+1)\geq s\}$ see \eqref{4lm2}, for $k_n= \bigg(\dfrac{p^*}{p}\bigg)^n-1$.

\end{lemma}
\begin{proof}

Remark that
\[v_m=0 \mbox{ where } u>0.\]
Then, proceeding analogously to Steps 1 and 2 in Lemma \ref{lemma2}, we end up with
\begin{align}
\nonumber
&\int_\Omega -b(u)|\nabla u|^{p-2}\nabla u \cdot \nabla v_m^{kp+1}-u v_m^{kp+1}
\\
\nonumber
&>\dfrac{a_0^{p-1}(kp+1)}{(k+1)^p}\|\nabla v_m^{k+1}\|^p_{L^p}+ \|v_m\|^{kp+2}_{L^{kp+2}}, \forall m>0 \mbox{ and } k\geq0.
\end{align}
Moreover,
\begin{align*}
\int_\Omega -f(x,u)v_m^{kp+1}&\leq C\bigg(\|u_-\|^{kp+1}_{L^{\alpha + kp+1}}+\|u_-\|^{\alpha+kp+1}_{L^{\alpha+kp+1}}\bigg)
\end{align*}
where $C=C(c_1, \alpha, \Omega)>0$.

Then, by combining the definition of weak solution of \eqref{p1}, the Sobolev Embedding Theorem and the Fatou Lemma, we have

\begin{align}
\nonumber
\|u_-\|_{L^{p^*(k+1)}}
&< \bigg[C\dfrac{(k+1)^p}{(kp+1)}\bigg(\|u_-\|^{kp+1}_{L^{\alpha+kp+1}}
\\
\nonumber
&+\|u_-\|^{\alpha+kp+1}_{L^{\alpha+kp+1}}+\|u_-\|^{(k+1)p}_{L^{(k+1)p}}\bigg)\bigg]^{1/p(k+1)}.
\end{align}

In addition, since $\alpha +kp+1<p(k+1)p$

\[\|u_-\|_{L^{\alpha+kp+1}}\leq \big(|\Omega|+1\big)\|u_-\|_{L^{k(p+1)}},\]
straightforward calculations lead us to
\begin{align}
\nonumber
\|u_-\|_{L^{p^*(k+1)}}&< \bigg[C\dfrac{(k+1)^p}{(kp+1)}\bigg(\|u_-\|^{kp+1}_{L^{(k+1)p}}
\\
\nonumber
&+\|u_-\|^{\alpha+kp+1}_{L^{(k+1)p}}+\|u_-\|^{(k+1)p}_{L^{(k+1)p}}\bigg)\bigg]^{1/p(k+1)}
\\
\nonumber
&\leq (k+1)^{1/(k+1)}\bigg[C\bigg(1+\|u_-\|^{(k+1)p}_{L^{(k+1)p}}\bigg)\bigg]^{1/p(k+1)}
\\
\label{2lm3}
&\leq e^{1/e}C^{1/p(k+1)}\bigg(1+\|u_-\|_{L^{(k+1)p}}\bigg),
\end{align}
where $C=C(N, p, a_0, \alpha, c_1, \Omega)>0$.


Now, remark that
\[ u_- \in L^s, \forall s \in (1,+\infty).\]

Indeed, since $k_n= \bigg(\dfrac{p^*}{p}\bigg)^n-1$, in particular,  $(k_1+1)p=p^*$. Then, by \eqref{2lm3},
\[\|u_-\|_{L^{p^*(k_1+1)}} \leq C(k_1)\bigg ( 1 + \|u_-\|_{L^{p^*}}\bigg).\]

Thus, since $(k_2+1)p=(k_1+1)p^*$, once again from inequality \eqref{2lm3} we arrive at
\[\|u_-\|_{L^{p^*(k_2+1)}} \leq C(k_2)\bigg ( 1 + \|u_-\|_{L^{p^*(k_1+1)}}\bigg) \leq C(k_2)\bigg ( 1 + \|u_-\|_{L^{p^*}}\bigg).\]

Hence, by induction, it follows that

\begin{equation}
\nonumber
\|u_-\|_{L^{p^*(k_n+1)}}\leq C(k_n)\bigg(1+\|u_-\|_{L^{p^*}}\bigg),
\end{equation}
what proves the claim.

Finally, observe that by the choice of $n_s$ and $k_0$, we have

\[
\|u_-\|_{L^{s}}\leq C\bigg(1+\|u_-\|_{L^{p^*}}\bigg)
\]
where $C=C(N, p, p^*, n_s, a_0, c_1, \alpha, \Omega)>0$.

Then, since it is clear that the latter inequalities also hold for $u_+$, with minor modifications on the latter arguments, the result follows.

\end{proof}

We are now ready to provide proofs for the core contributions of the present section. Actually, the next results guarantee $L^\infty$ and fractional order a priori bounds for solutions of \eqref{p1}. Once again, we begin with the case where $p<N$ and $p-1\leq \alpha <p^*-1$.

\begin{proposition} \label{prop2} Under hypotheses \eqref{H1} and \eqref{H2}, suppose that, $2<p<N$ and $u\in W^{1,p}$ is a weak solution of \eqref{p1}. If  $R\in \mathbb{R}$ is a fixed number satisfying
\begin{equation}
\nonumber
\label{exponents1}
 R>\max\{2, (N-4)/(p-2)\},
\end{equation}
and
\[ p-1\leq\alpha<p^*-1\]
set
\[q=\alpha R.\]
Then, $u\in L^\infty$ and also
\begin{align}
\label{1pr2}
\|u\|_{L^\infty} &\leq C \bigg(1+a(\|u\|^p_{W^{1,p}})^{\omega_{p,R}}+\big(h_q(\|u\|_{L^{p^*}})\big)^\alpha
\\
\nonumber
&+a(\|u\|^p_{W^{1,p}})^{\omega_{p,R}}\big(h_q(\|u\|_{L^{p^*}})\big)^{q/r_R}\bigg),
\end{align}
where $C=C(N, p, p^*, q, R, a_0,  c_1, \alpha, \delta, \Omega)>0$, $r_R$ is given in \eqref{exponents2}, $\omega_{p,r_R}$ in \eqref{exponents4}, $h_q$ is defined in \eqref{2lm2}, and $\delta$ is the same as in Lemma \ref{lemma2}.

Moreover, $u \in \mathcal{N}^{1+\frac{2}{r_s},r_s}$ and
\begin{equation}
\label{2pr2}
\|u\|_{\mathcal{N}^{1+\frac{2}{r_s},r_s}} \leq C \bigg(1+a(\|u\|^p_{W^{1,p}})^{T_s}+\big(1+a(\|u\|^p_{W^{1,p}})\big)^{T_s}h_q(\|u\|_{L^{p^*}})^{R_s}\bigg),
\end{equation}
for every $s \in (2,+\infty)$, where
\[R_s=\max \bigg\{\alpha^2,\dfrac{\alpha q}{ r_R},\dfrac{\alpha^2 s}{r_s}, \dfrac{\alpha q s}{r_R r_s}\bigg\}, T_s=\omega_{p,s}+\dfrac{\alpha \omega_{p,R}(r_s+s)}{r_s},\]
and
$C=C(N, p, p^*, q, R, R_s, T_s, a_0, c_1,  \alpha, \delta, \Omega)>0.$

\end{proposition}
\begin{proof}

For the sake of clarity, we split this proof into two steps.

\noindent \textit{Step 1.} {\it $L^\infty$ estimates}.

Owing to Lemma \ref{lemma2}, we readily see that  $u \in L^q$ and $|u|^\alpha \in L^R$.

Then, set $\rho=r_{R}$, and observe that by combining \eqref{H2} and Proposition \ref{prop1}, we conclude that $u\in \mathcal{N}^{1+\frac{2}{\rho},\rho}$ and
\begin{equation}
\nonumber
\|u\|_{\mathcal{N}^{1+\frac{2}{\rho},\rho}}
\leq C\bigg( 1+\|u \|^\alpha_{L^q}+a(\|u\|^p_{W^{1,p}})^{\omega_{p,R}}\|u\|_{L^q}^{q/\rho}\bigg),
\end{equation}
since $p>2$. In addition, by inequality \eqref{1lm2}, p. \pageref{1lm2}, we have that
\begin{align}
\nonumber
\|u \|^\alpha_{L^q}+a(\|u\|^p_{W^{1,p}})^{\omega_{p,R}}\|u\|_{L^q}^{q/\rho} \leq& C \bigg(1+a(\|u\|^p_{W^{1,p}})^{\omega_{p,R}}+h_q(\|u\|_{L^{p^*}})^\alpha
\\
\nonumber
&
 +a(\|u\|^p_{W^{1,p}})^{\omega_{p,R}}h_q(\|u\|_{L^{p^*}})^{q/\rho}\bigg)
\\
\nonumber
\end{align}
where \[h_q=x+x^{(n_q+1)\lambda},x\geq 0\]
and $n_q$ is the same as in Lemma \ref{lemma2}, see \eqref{4lm2} and \eqref{sequence1}, p. \pageref{sequence1}.

Then, by plugging the latter inequalities we get
\begin{align}
\label{3pr2}
\|u\|_{\mathcal{N}^{1+\frac{2}{\rho},\rho}}&\leq C\bigg(1+a(\|u\|^p_{W^{1,p}})^{\omega_{p,R}}+h_q(\|u\|_{L^{p^*}})^\alpha
\\
&+a(\|u\|^p_{W^{1,p}})^{\omega_{p,R}}h_q(\|u\|_{L^{p^*}})^{q/\rho}\bigg),
\nonumber
\end{align}
for $C=C(N,p,R,\Omega)>0$.

Nevertheless, recall that for every $\epsilon>0$, sufficiently small, there holds
\[
\mathcal{N}^{1+\frac{2}{\rho},\rho} \hookrightarrow W^{1+\frac{2}{\rho}-\epsilon,\rho} \hookrightarrow C^{0,\sigma_\epsilon}(\overline{\Omega}),
\]
where
\[\sigma_\epsilon = 1 +\dfrac{2-N}{\rho}-\epsilon <1+\dfrac{2-N}{\rho}.\]
For instance, see  \cite{knees} Lemma 2.1, \cite{grisvard} Thm. 1.4.4.1 and the subsequent commentaries. In addition, observe that by the choice of $R$ and $\rho$
 \[
 1+\dfrac{2-N}{\rho}>0,
 \]
 so that, there exists $\epsilon>0$ for which $\sigma_\epsilon>0$.

Further, by fixing $\epsilon=\frac{\rho+(2-N)}{2\rho}$, from \eqref{3pr2}, there follows that

\begin{align}
\label{4pr2}
\|u\|_{L^\infty}&\leq C\|u\|_{C^{0,\frac{1}\sigma_{\rho}}(\overline{\Omega})}
\\
\nonumber
&\leq C\bigg(1+a(\|u\|^p_{W^{1,p}})^{\omega_{p,R}}+h_q(\|u\|_{L^{p^*}})^\alpha
\\
\nonumber
&+a(\|u\|_{W^{1,p}}^p)^{\omega_{p,R}}h_q(\|u\|_{L^{p^*}})^{q/\rho}\bigg).
\end{align}

\noindent \textit{Step 2.} {\it Fractional order estimates}.

Now, given $s$, such that $2<s<+\infty$, recall that $r_s=s(p-2)+2$, see \eqref{exponents2}. Then, by combining Proposition \ref{prop1} and inequality \eqref{1pr2}, since $|u|^\alpha \in L^s$ and $p>2$,  we get
\begin{align}
\nonumber
\|u\|_{\mathcal{N}^{1+\frac{2}{r_s},r_s}} 
&\leq C\bigg(1+ \|u\|^\alpha_{L^{\infty}}+a(\|u\|^p_{W^{1,p}})^{\omega_{p,s}}\|u\|^{\alpha s/r_s}_{L^{\infty}} \bigg)
\\
\nonumber
&\leq C\bigg(1+a(\|u\|^p_{W^{1,p}})^{\alpha \omega_{p,R}}+h_q(\|u\|_{L^{p^*}})^{\alpha^2}
\\
\nonumber
&+a(\|u\|^p_{W^{1,p}})^{\alpha \omega_{p,R}}h_q(\|u\|_{L^{p^*}})^{\alpha q/\rho}+a(\|u\|^p_{W^{1,p}})^{\omega_{p,s}}
\\
\nonumber
&+a(\|u\|^p_{W^{1,p}})^{\omega_{p,s}+(\alpha s  \omega_{p,R})/r_s}+a(\|u\|^p_{W^{1,p}})^{\omega_{p,s}}h_q(\|u\|_{L^{p^*}})^{\alpha^2 s/r_s}
\\
\nonumber
&
+a(\|u\|^p_{W^{1,p}})^{\omega_{p,s}+(\alpha s  \omega_{p,R})/r_s}\big(h_q(\|u\|_{L^{p^*}})\big)^{\alpha qs/\rho r_s}\bigg)
\\
\nonumber
&\leq C \bigg(1+a(\|u\|^p_{W^{1,p}})^{T_s}+\big(1+a(\|u\|^p_{W^{1,p}}\big)^{T_s}h_q(\|u\|_{L^{p^*}})^{R_s}\bigg)
\end{align}
what guarantees the validity of \eqref{2pr2}, completing this proof.
\end{proof}

In the previous results,  $L^\infty$ and fractional order estimates for solutions of \eqref{p1} were given in the case when $p-1\leq \alpha < p^*-1$, in terms of $h(\|u\|_{L^{p^*}})$. In the case where $1\leq \alpha < p-1$, naturally, these results can be  improved.

\begin{proposition} \label{prop3} Under hypotheses \eqref{H1} and \eqref{H2}, suppose that $2<p<N$, $u\in W^{1,p}$ is a weak solution of \eqref{p1}, and that $1\leq \alpha < p-1$. Further, consider
\[q=\alpha R\]
for $R$ as in \eqref{exponents1}.

Then, there holds that

\begin{equation}
\label{1lm3}
\|u\|_{L^\infty} \leq C\bigg( 1+a(\|u\|^p_{W^{1,p}})^{\omega_{p,R}}+\|u\|^\alpha_{L^{p^*}}+a(\|u\|^p_{W^{1,p}})^{\omega_{p,R}}\|u\|^{q/r_R}_{L^{p^*}}\bigg),
\end{equation}

where $r_R$, is given in \eqref{exponents2}, and $C=C(N, p, p^*, R, a_0, c_1, \alpha, \Omega)>0$.

Moreover, given $s\in (2,+\infty)$, $u \in \mathcal{N}^{1+\frac{2}{r_s},r_s}$ and
\begin{align}
\nonumber
\|u\|_{\mathcal{N}^{1+\frac{2}{r_s},r_s}} &\leq C\bigg(1+\|u\|_{L^{p^*}}^{R_s}\bigg),
\end{align}
for \[  R_s=\max \bigg\{\alpha^2,\dfrac{\alpha q}{ r_R},\dfrac{\alpha^2 s}{r_s}, \dfrac{\alpha q s}{r_R r_s}\bigg\}\mbox{ and }T_s=\omega_{p,s}+\dfrac{\alpha \omega_{p,R}(r_s+s)}{r_s},\]
where
$C=C(N, p, p^*, R, R_s, T_s, a_0, c_1, \alpha, \Omega)>0$.
\end{proposition}
\begin{proof}
The argument is analogous to the proof of Lemma \ref{lemma2}. For the convenience of the reader, we exhibit its details.


We now consider
\[
v_m=\min\{u_-,m\},
\]
where $v_m \in L^\infty\cap W^{1,p}$ for all $m>0$. Given $k\geq0$, we are going to use $-v_m^{kp+1}$ as a test function in \eqref{p1}.

\noindent \textit{Step 1.} {\it Basic estimates}.
Then, recalling that $r_R=R(p-2)+2$, once again  by combining \eqref{H2} and Proposition \ref{prop1}, we conclude that $u\in \mathcal{N}^{1+\frac{2}{r_R},r_R}$ and

\begin{align}
\nonumber
\|u\|_{\mathcal{N}^{1+\frac{2}{r_R},r_R}}
&\leq C\bigg( 1+\|u \|^\alpha_{L^q}+a(\|u\|^p_{W^{1,p}})^{\omega_{p,R}}\|u\|_{L^q}^{q/r_R}\bigg),
\end{align}
since $p>2$.

Hence, due to an argument analogous to \eqref{4pr2}, by recalling Lemma \ref{lemma4}, we prove \eqref{1lm3}.


\noindent \textit{Step 2.} {\it Fractional order estimates}.

For this final step, given $s$, such that $2<s<+\infty$, by Proposition \ref{prop1}, as $|u|^\alpha \in L^s$, we get
\begin{align}
\nonumber
\|u\|_{\mathcal{N}^{1+\frac{2}{r_s},r_s}} &\leq C\bigg( 1+\|u\|^\alpha_{L^{\alpha s}}+a(\|u\|^p_{W^{1,p}})^{\omega_{p,s}}\|u\|^{\alpha s/r_s}_{L^{\alpha s}}\bigg)
\\
\nonumber
&\leq C\bigg(1+ \|u\|^\alpha_{L^{\infty}}+a(\|u\|^p_{W^{1,p}})^{\omega_{p,s}}\|u\|^{\alpha s/r_s}_{L^{\infty}} \bigg)
\\
\nonumber
&\leq C\bigg(1+a(\|u\|^p_{W^{1,p}})^{\alpha\omega_{p,R}}+\|u\|^{\alpha^2}_{L^{p^*}}+a(\|u\|^p_{W^{1,p}})^{\alpha\omega_{p,R}}\|u\|^{\alpha q/r_R}_{L^{p^*}}
\\
\nonumber
&+a(\|u\|^p_{W^{1,p}})^{\omega_{p,s}+\alpha\omega_{p,R}s/r_s}+a(\|u\|^p_{W^{1,p}})^{\omega_{p,s}}\|u\|^{\alpha^2 s/r_s}_{L^{p^*}}
\\
\nonumber
&+a(\|u\|^p_{W^{1,p}})^{\omega_{p,s}+\alpha\omega_{p,R}s/r_s}\|u\|_{L^{p^*}}^{\alpha qs/r_Rr_s}\bigg)
\\
\nonumber
&\leq C \bigg(1+a(\|u\|^p_{W^{1,p}})^{T_s}+\big(1+a(\|u\|^p_{W^{1,p}})\big)^{T_s}\|u\|_{L^{p^*}}^{R_s}\bigg),
\end{align}
and the result follows.
\end{proof}

 By means of minor modifications on the arguments,  we can also relax hypothesis (H2) and still, guarantee the validity of $L^\infty$ and fractional order estimates for the cases where $p\geq N$. Indeed, for $p \geq N$, we suppose that \eqref{H2P} holds, and by straightforward modifications on the proofs of the latter propositions, we obtain the following result.

\begin{proposition} \label{prop4} Under hypotheses \eqref{H1} and \eqref{H2P}, consider  $u\in W^{1,p}$, a weak solution of \eqref{p1}.

Then, there holds that,  if $p=N$,

\begin{equation}
\label{1prop4}
\|u\|_{L^\infty} \leq C\bigg( 1+\|u\|^\alpha_{W^{1,p}}+\|u\|^{q/r}_{W^{1,p}}\bigg),
\end{equation}

for $r_{R}$ defined in \eqref{exponents2}, $\omega_{p,R}$ in \eqref{exponents4}, and $C=C(N, p, R, a_0, n_0, \alpha, c_1, \beta, \Omega)>0$.

Moreover, for every $p\geq N$, given $s\in (2,+\infty)$ there holds that
\begin{align}
\label{2prop4}
\|u\|_{\mathcal{N}^{1+\frac{2}{r_s},r_s}} &\leq C\bigg(1+a(\|u\|_{W^{1,p}}^p)^{\omega_{p,s}}+\|u\|_{W^{1,p}}^{\alpha}
\\
\nonumber
&+a(\|u\|_{W^{1,p}}^p)^{\omega_{p,s}}\|u\|_{W^{1,p}}^{\alpha s/r_s}\bigg) \mbox{ if } p=N,
\\
\nonumber
\|u\|_{\mathcal{N}^{1+\frac{2}{r_s},r_s}} &\leq C\bigg(1+\|u\|_{W^{1,p}}^{\alpha}+a(\|u\|_{W^{1,p}}^p)^{\omega_{p,s}}\|u\|_{W^{1,p}}^{\alpha s/r_s}\bigg) \mbox{ if } p>N,
\end{align}
for $C=C(N, p, R, a_0, n_0,\alpha, c_1, \Omega)>0$.
\end{proposition}
\begin{proof}
For \eqref{1prop4}, we only have to mimic the argument used in {\it Step. 2} on the proof of Proposition \ref{prop3}, replacing $L^{p^*}$ norms by $W^{1,p}$ norms.
In an analogous manner, by means of subtle modifications on the proof of {\it Step. 3} in Proposition \ref{prop3}, and by using the Sobolev Embedding Theorem, we obtain \eqref{2prop4}.

\end{proof}

Now, we are in the position to prove the aforementioned uniform $L^\infty$ a priori bounds for the solutions of \eqref{p1}.
We once more stress that this is where hypotheses \eqref{NQ1} and \eqref{H3} take place.

\begin{proposition} \label{prop5} Suppose that $p>2$. Under hypotheses \eqref{H1}, \eqref{H2} or \eqref{H2P}, and \eqref{H3}, suppose that \eqref{NQ1} holds.
Then, if
 \begin{equation}\label{sigma=beta2}
	 \nu c_3 < c_2 \mbox{ and } \sigma \geq \beta;
 \end{equation}
there exists $C>0$, for which
\begin{equation}
\label{1prop5}
\|u\|_{L^\infty} \leq C,
\end{equation}

for every weak solution of \eqref{p1}, where $C(N, p, \|a\|_{L^\infty_{loc}},\|F\|_{L^\infty_{loc}}, \alpha, \beta, \nu, \sigma, \Omega)>0$.

\end{proposition}
\begin{proof}
By Propositions \ref{prop2}, \ref{prop3} and \ref{prop4}, it is enough to prove that
\[\|u\|_{W^{1,p}} \leq C,\]
for every $u$ weak solution of \eqref{p1}.
Since \eqref{NQ1} holds, given $\epsilon=\epsilon(c_3)>0$, there exists
\[ R_1=R_1(f,\sigma,c_2)>0\]
such that
\[ f(x,t)t-\nu F(x,t) \leq (-c_2+\epsilon)|t|^\sigma \mbox{ a.e in } \Omega\]
for all $|t|\geq R_1$.

In addition, given $\delta=\delta(c_3)>0$, by \eqref{H3}, there exists
\[R_2=R_2(\beta,f,c_3)>0\]
such that
\[ F(x,t) \leq (c_3+\delta)|t|^\beta \mbox{ a.e. in } \Omega\]
for all $|t|\geq R_2$.

Hence, by the latter inequalities
\begin{equation}
\label{2prop5}
f(x,t)t \leq (-c_2+\epsilon)|t|^\sigma +\nu(c_3+\delta)|t|^\beta\mbox{ a.e in } \Omega.
\end{equation}

Now we are going to choose carefully the constants $\epsilon$ and $\delta$.
Choose $\delta$ and $\epsilon$ such that $\epsilon + \nu \delta < c_2-\nu c_3$. Thence, observe that if \eqref{sigma=beta2} holds, clearly there exists
\[
R_3=R_3(c_2,c_3,\beta,\nu, \sigma)>0
\]
such that
\begin{equation}
\label{3prop5}
-(c_2-\epsilon)|t|^\sigma+\nu(c_3+\delta)|t|^\beta\leq 0
\end{equation}
for all $|t|>R_3$.

Thus, set
\[R=\max_{1\leq i\leq 3} \{R_i,1\}\]
and
\[ c= \max_{x\in \Omega, t\leq |R|} |f(x,t)|,\]
which is finite by \eqref{H2}.

Moreover, given $u\in W^{1,p}$ denote

\[ \Omega_-=\{x:|u(x)|\leq R\} \mbox{ and } \Omega_+=\{x:|u(x)|> R\}.\]

Hence, by the choice of $R$, $c$, combining \eqref{2prop5} and \eqref{3prop5}, we arrive at

\begin{align*}
\int_\Omega f(x,u)u &= \int_{\Omega_+} f(x,u)u+\int_{\Omega_-} f(x,u)u
\\
&\leq  \int_{\Omega_+} -(c_2-\epsilon)|t|^\sigma +\nu(c_3+\delta)|t|^\beta+\int_{\Omega_-} |f(x,u)|u|
\\
&\leq  \int_{\Omega_-} cR =K
\end{align*}

Thus, \eqref{1prop5} holds for every weak solution of \eqref{p1}.
\end{proof}

In an analogous manner, we shall employ hypotheses \eqref{NQ1} and \eqref{H3P} in order to obtain $L^\infty$ uniform bounds.
\begin{proposition} \label{prop6} Suppose that $p>2$. Under hypotheses \eqref{H1}, \eqref{H2} or \eqref{H2P}, and \eqref{H3P}, suppose that \eqref{NQ1} holds.
Then, we obtain the same conclusion of Proposition \ref{prop5}.
\end{proposition}
\begin{proof}
By Propositions \ref{prop2}, \ref{prop3} and \ref{prop4}, it is enough to prove that
\[\|u\|_{W^{1,p}} \leq C,\]
for every $u$ weak solution of \eqref{p1}.
Since \eqref{NQ1} holds, given $\epsilon=\epsilon(c_3)>0$, there exists
\[ R_1=R_1(f,\sigma,c_2)>0\]
such that
\[ f(x,t)t-\nu F(x,t) \leq (-c_2+\epsilon)|t|^\sigma \mbox{ a.e in } \Omega\]
for all $|t|\geq R_1$.

In addition, given $\delta>0$, by \eqref{H3P}, there exists
\[R_2=R_2(\beta,f,\delta)>0\]
such that
\[ -F(x,t) \geq \delta|t|^\beta \mbox{ a.e. in } \Omega\]
for all $|t|\geq R_2$.

Hence, by the latter inequalities
\begin{equation}
\nonumber
f(x,t)t \leq (-c_2+\epsilon)|t|^\sigma -\nu\delta|t|^\beta\mbox{ a.e in } \Omega.
\end{equation}
Choose $\epsilon$ such that $-c_2+\epsilon <0$. Thence, clearly
\begin{equation}
\nonumber
(-c_2+\epsilon)|t|^\sigma -\nu\delta|t|^\beta\leq 0
\end{equation}
for all $|t|>R_2$. The conclusion is the same as in Proposition \ref{prop5}.
\end{proof}

Now, with the tools of the last section in hand, we are able to prove our main results.

\section{Proof of the main Results}\label{proof}

\textbf{Proof of Theorem \ref{thm1}.}
Once again, the proof of the existence of solution will be based on the Leray--Schauder Fixed Point Theorem.

Beforehand, let us choose $\tau>2$ such that for, there holds $r_\tau^*/\tau>\alpha$.
In this fashion, remark that by the choice of $\tau$, given $v \in W^{1,r_\tau}$, there holds that $f(x,v)\in L^\tau$.

Now, we consider the following family of operators
\[T_t:W^{1,r_\tau}\to W^{1,r_\tau} \mbox{ for } 0\leq t \leq 1, \mbox{ where } T_t(v)=u \mbox{ if and only if }\]

\begin{equation}
\label{1thm1}
\left\{
\begin{aligned}
-b(v)\Delta_p u + u & =  t f(x,v) && \text{in } \Omega,\\
     \dfrac{\partial u}{\partial\eta} & =    0 && \text{on } \partial \Omega.
\end{aligned}
\right.
\end{equation}

which are well--defined by Proposition \ref{prop1}, since there exists a unique \[u\in \mathcal{N}^{1+\frac{2}{r_\tau}, r_\tau}\] satisfying \eqref{1thm1}.

 At this point, let us stress that by combining Steps 1,2,3 and 5 of aforementioned proposition with Propositions \ref{prop2}-\ref{prop4}, it is straightforward to prove that
$T_t(.)$ is continuous, compact and uniformly continuous with respect to $t$ on bounded sets of $W^{1,r_\tau}$, what guarantees the existence of a solution $u$.
However, for the reader convenience, we will check the validity of the hypotheses of Leray--Schauder's theorem in four steps.

\subsection*{ Step 1. \textit{$T_t$ is compact.}} It is enough to remark that for a given $\{v_n\}_{\mathbb{N}} \subset W^{1,r_\tau}$, bounded, for $u_n=T_t(v_n)$ we obtain by Propositions \ref{prop2}-\ref{prop4} that $\{u_n\}_{\mathbb{N}}$ is bounded in $\mathcal{N}^{1+\frac{2}{r_\tau},r_\tau}$ so that $T_t$ is obviously compact.

\subsection*{ Step 2. \textit{$T_t$ is continuous.}} Since if $v_n \to v$ in $W^{1,r_\tau}$,  it is clear that
\[\lim_{n\to+\infty} \int_\Omega f(x,v_n) \phi = \int_\Omega f(x,v)\phi, \ \forall \phi\in W^{1,p}\]
then, by repeating the same argument given in Step 2 of Proposition \ref{prop1}, we prove our claim.

\subsection*{ Step 3. \textit{$T_t$ is uniformly continuous with respect to $t$ on bounded subsets of $W^{1,r_\tau}$.}} Since $r_\tau>p$, by Propositions \ref{prop2}-\ref{prop4} there exists $C=C(A)>0$ for which
\[ \|u\|_{\mathcal{N}^{1+\frac{2}{r_\tau}, r_\tau}}\leq C\]
for every $u=T_t(v)$ where $v \in A \subset W^{1,r_\tau}$, a bounded subset.

Further, by mimicking the proof of Step 3, we obtain that, for $u_i=T_{t_i} v$, where $v \in A $, we have
\[\|u_1-u_2\|_\wpp \leq C|t_1-t_2|^{1/p}\|f(x,v)\|_{L^\tau}^{2/p} + C|t_1-t_2|^{1/2}\|f(x,v)\|_{L^\tau}.\]

In addition, by recalling that for $\epsilon>0$ sufficiently small,
\[T_t(A) \subset \mathcal{N}^{1+\frac{2}{r_\tau}, r_\tau} \hookrightarrow\hookrightarrow W^{1,r_\tau+\epsilon},\]
and by combining the latter inequalities with the Gagliardo--Nirenberg inequality,  we obtain that
\begin{align*}
\|u_1-u_2\|_{W^{1,r_\tau}} &\leq C\|u_1-u_2\|^\gamma_{W^{1,p}}\|u_1-u_2\|^{1-\gamma}_{W^{1,r_\tau+\epsilon}}
\\
&\leq C\bigg(|t_1-t_2|^{\gamma/p}+|t_1-t_2|^{\gamma/2}\bigg),
\end{align*}
what guarantees that $T_t$ is uniformly continuous with respect to $t$ on bounded subsets of $W^{1,r_\tau}$.

It remains to prove that the fixed--points of $T_t$ are uniformly bounded in $W^{1,r_\tau}$, what we address in the next step.

\subsection*{ Step 4. \textit{A priori bounds and uniqueness}.} We claim that there exists $C>0$ such that for every fixed point of $T_1$
\[\|u\|_{W^{1,r_\tau}} \leq C.\]

In fact, by Proposition \ref{prop5} with $p=r_\tau$, there exists $C>0$ such that for all $u\in W^{1,r_\tau}$ solution of $T_t (u)=u$ there holds that $\|u\|_{W^{1,r_\tau}}\leq C$. Therefore, by Propositions \ref{prop2}, \ref{prop3} and \ref{prop4}, whereas $a(.)$ is a continuous functional, we arrive at
\[\|u\|_{\mathcal{N}^{1+\frac{2}{r_\tau}, r_\tau}} \leq C.\]

 Further, since it is clear that $T_0(v)\equiv 0$, that is, we have uniqueness of solutions when $t=0$, as a consequence of the Leray--Schauder fixed point theorem, there exists $u \in W^{1, r_\tau}$, a weak solution of \eqref{p1}.

For the nontriviality of this solution, remark that once $u=T_1(u)$ is given by Proposition \ref{prop1}, then $u$ satisfies \eqref{p1} a.e. in $\Omega$.
Complementarily, since for every $c\in \mathbb{R}$
\[\mu\bigg(\{x\in \Omega: f(.,c)\neq c\}\bigg)>0,\] we have that $u\not\equiv c $ in $\Omega$ for all $c\in \mathbb{R}$, so that $u$ is a nontrivial strong solution of \eqref{p1}.

Finally, by the construction of $T_t$, and by Proposition \ref{prop5}, inequality \eqref{1prop5}, there exists $C>0$ for which
\[ \|u\|_{L^\infty}\leq C \mbox{ and } \|u\|_{W^{1,r_s}}\leq C.\]
Therefore, by combining the latter inequalities and Propositions \ref{prop2}-\ref{prop4}, for every $s\in (2,+\infty)$, there holds
\[\|u\|_{\mathcal{N}^{1+\frac{2}{r_s},r_s}}\leq C,\]
what completes this proof.

\subsection*{Proof of Theorem \ref{thm2}}
Remark that this proof is completely analogous to the previous one. In fact, we will employ the Leray--Schauder Fixed Point Theorem one more time, for the same family of operators $T_t(.)$.

In this fashion,  we  need to prove the validity of the a priori bounds, i.e., that there exists $C>0$ which for every $u$ fixed point of $T_1$
\[\|u\|_{W^{1,r_\tau}} \leq C.\]
In turn, the arguments is the same as in Theorem \ref{thm1}, however, instead of Proposition \ref{prop5}, we employ Proposition \ref{prop6}.

 Thus, since the operator is the same as in Theorem \ref{thm1}, by repeating the proof that the $T_t(.)$ is continuous, compact and uniformly continuous with respect to $t$, we obtain the existence of a solution.

At last, observe that the proof that the solution is nontrivial is exactly the same, since we assume the same hypotheses of nontriviality (\ref{H2}) and (\ref{H2P}).

\section{Appendix}\label{appendix}
The purpose of the present section is to exhibit some of the basic background on Nikolskii spaces and fractional regularity, specially, concerning the aspects of the theory which are somehow connected to the investigation of solutions to the $p$-Laplacian.
For the interested reader, without the intention of being complete, we recommend the excellent monographs \cite{bin,kufner,leoni,leoni2,rsc,saw,triebel3}, and the references thereof, which cover the subject discussed in this section in detail.

\subsection{ Definition of $\mathcal{N}^{\sigma,r}$ and relations with Besov Spaces }
At this point, for the sake of clarity, we recall the basic framework of fractional spaces which appear in this work. First, we state our definition of Nikolskii spaces.

\begin{definition}
\label{nikolskii}
Consider $\Omega \subset \mathbb{R}^N$, a bounded smooth domain, $\sigma \in (1,2)$ and $r> 1$. We define the Nikolskii space $\mathcal{N}^{\sigma,r}(\Omega)$ by
\[\mathcal{N}^{\sigma,r}(\Omega)=\bigg\{u\in W^{1,r}(\Omega): \sum_{i=1}^N\sup_{h\neq 0} \bigg{(}\int_{\Omega_{|h|}}\dfrac{\big{|}\partial_{x_i} u(x+h)-\partial_{x_i} u(x)\big{|}^r}{|h|^{\sigma r}}\bigg{)}^{1/r}<+\infty\bigg\}\]
where $\Omega_{|h|}=\{x\in \Omega: dist(x,\partial \Omega)<|h|\}$ and $h\in\mathbb{R}^N$.
Moreover, consider the following norm
\[\|u\|_{\mathcal{N}^{\sigma,r}}=\|u\|_{W^{1,r}}+[[u]]_{\mathcal{N}^{\sigma,r}},\]
where
\[[[u]]_{\mathcal{N}^{\sigma,r}}=\sum_{i=1}^N\sup_{h\neq 0} \bigg{(}\int_{\Omega_{|h|}}\dfrac{\big{|}\partial_{x_i} u(x+h)-\partial_{x_i} u(x)\big{|}^r}{|h|^{\sigma r}}\bigg{)}^{1/r}\]
is the so--called Gagliardo--Nikolskii seminorm.
\end{definition}
Remark that by the equivalence of norms in $\mathbb{R}^N$, $[[u]]_{\mathcal{N}^{\sigma,r}}$ is equivalent to

\[\sup_{h\neq 0} \bigg{(}\int_{\Omega_{|h|}}\dfrac{\big{|}\nabla u(x+h)-\nabla u(x)\big{|}^r}{|h|^{\sigma r}}\bigg{)}^{1/r}.\]

We stress that, from now on, for the rest of this section, we are going to consider $\Omega \subset \mathbb{R}^N$, a bounded smooth domain, $\sigma \in (1,2)$ and $r> 1$. Further, for the sake of simplicity we are going to denote $\mathcal{N}^{\sigma,r}(\Omega)$ simply by $\mathcal{N}^{\sigma,r}$ .

There are other characterizations of these spaces, we have chosen the latter one because it suits better to the context of solutions to degenerate equations with $p$-structure. The reader is invited to look at \cite{bin}, Section 18, \cite{rsc}, Sections 2.3 and 2.4, \cite{triebel3}, Section 1.1, \cite{leoni2}, Chapter 17, Sections 17.2, 17.6 -- specially Corollary 17.68 and Theorem 17.69 -- and 17.7, or in \cite{leoni}, Chapter 14.  Moreover, we recall that Nikolskii spaces are a particular case of the so--called Besov spaces, $B^{\sigma}_{r,q}$, which also generalize the Sobolev--Slobodeckii spaces, $W^{\sigma,r}$. For instance, under the conditions of Definition \ref{nikolskii},  there holds that
\[B^{\sigma}_{r,\infty}=\mathcal{N}^{\sigma,r} \mbox{ and } B^{\sigma}_{r,r}=W^{\sigma,r},\]
see \cite{bin}, p. 59, \cite{rsc}, Section 2.1, pp. 11--14, or \cite{leoni2}, pp. 539--540, and also \cite{leoni}, Section 14.8.
\subsection{ Basic properties}
Now we sketch the basic properties of these spaces. The proofs of the standard results are going to be left to the reader.

As one should expect, such spaces are complete with respect to the latter norm, for instance see \cite{bin}, Theorem 18.3, and \cite{leoni2}, p. 542, or Proposition 14.3 in \cite{leoni}.
\begin{lemma}
\label{banach}
The space $\mathcal{N}^{\sigma,r}$ endowed with the norm $\|.\|_{\mathcal{N}^{\sigma,r}}$ is a Banach space.
\end{lemma}

The next lemma, guarantees a compactness principle which is essential in the context of the fractional regularity theory used in the text.
\begin{lemma}
\label{compactness}
Suppose that $1<\sigma<2$ and $r>1$.
Then, given $\epsilon \in (0,\sigma]$ the embeddings
\[ \mathcal{N}^{\sigma,r}\hookrightarrow\hookrightarrow W^{\sigma-\epsilon,r} \hookrightarrow \mathcal{N}^{\sigma-\epsilon,r},\]
hold true.

In particular, there holds that
\[
\mathcal{N}^{\sigma,r} \hookrightarrow\hookrightarrow W^{1,r} \]
is a compact embedding.
\end{lemma}
For these embeddings we refer the interested reader to \cite{rsc}, Theorems 1 and 2, p. 82, to \cite{leoni2}, pp. 559--561, and also to \cite{knees}, Lemma 2.1.

\subsection{Fractional Regularity Tools and Example}
\label{toolsexample}
Now we address a fundamental nonlinear estimate responsible for linking the natural energy estimates appearing in $p$-Laplacian--like equations with Nikolskii spaces.

\begin{lemma}\label{normalemma}  Consider $r>2$ and $u\in W^{1,p}$ such that the Hessian $D^2u$ exists a.e. in $\Omega$ and
\begin{equation}
\nonumber
\int_{\Omega} |\nabla u|^{r-2}|D^2 u|^2<+\infty.
\end{equation}
Then,  there exists $C=C(\Omega,N,p)>0$ for which
\begin{equation*}
 [[u]]^r_{\mathcal{N}^{1+\frac{2}{r},r}}\leq C\int_{\Omega} |\nabla u|^{r-2}|D^2 u|^2.
\end{equation*}
\end{lemma}

\begin{proof}

As a first step, remark that  $|\nabla u|^{(r-2)/2}\nabla u \in W^{1,2}$. In fact,
\begin{equation}
\nonumber
\big{|}|\nabla u|^{(r-2)/2}\nabla u\big{|}^2\leq |\nabla u|^p \in L^1,
\end{equation}
and,  for $i=1,..., N, $  a.e. in $\Omega$, we have
\begin{eqnarray}
\nonumber
\bigg{|}\dfrac{\partial}{\partial x_i}\bigg{(}|\nabla u|^{(r-2)/2}\dfrac{\partial u}{\partial x_j}\bigg{)}\bigg{|} &\leq& |\nabla u|^{(r-2)/2}\bigg{|}\dfrac{\partial^2 u}{\partial x_i\partial x_j}\bigg{|} \nonumber \\ &+& \dfrac{r-2}{2}|\nabla u|^{(r-6)/2} \sum_{k=1}^N |\nabla u|^2\bigg{|}\dfrac{\partial^2 u}{\partial x_k \partial x_i}\bigg{|}
\nonumber \\
\label{2lemma1}
&\leq& C |\nabla u|^{(r-2)/2}|D^2 u| \in L^2.
\end{eqnarray}

Now, we proceed to the most important nonlinear estimate of this proof. Actually, we claim that
\begin{equation}
\label{5lemma1}
[[u]]_{\mathcal{N}^{1+2/r,r}}^r \leq C  \sup_{h > 0 } \int_{\Omega_{|h|}}\big{|}D^h (| \nabla u|^{(r-2)/2} \nabla u)\big{|}^2,
\end{equation}
 where $$D^h g= \frac{g(.+h)-g(.)}{|h|}$$ is the difference quotient of $g$ in the direction $h \in \mathbb{R}^n$, $ h \neq 0$.

 Indeed, recall  that for all $U, \ V \in \mathbb{R}^N$,  there exists $C>0$ such that
\begin{equation}
\nonumber
|U-V|^r\leq C \bigl{|}U|U|^{(r-2)/2}-V|V|^{(r-2)/2}\bigr{|}^2.
\end{equation}

Thence, let us employ the latter inequality to  $$U= T^h \nabla u(.) =\nabla u(.+h) \ \text{ and for } \  V=\nabla u(.),$$ so that
\begin{eqnarray*}
|T^h \nabla u-\nabla u|^r&\leq &
C \ \big{|}T^h \nabla u|T^h \nabla u|^{(r-2)/2}-\nabla u|\nabla u|^{(r-2)/2}\big{|}^2
\\
\nonumber
&= & C \ |D^h(|\nabla u|^{(r-2)/2}\nabla u)|^2 |h|^2.
\end{eqnarray*}
Hence, by dividing by $|h|^2$, after an integration over $\Omega_{|h|}$, we arrive at
\begin{equation*}
\int_{\Omega_{|h|}}\dfrac{|T^h \nabla u-\nabla u|^r}{|h|^2}\leq C \int_{\Omega_{|h|}}\big{|}D^h (| \nabla u|^{(r-2)/2} \nabla u)\big{|}^2,
\end{equation*}
which implies \eqref{5lemma1}.

Further, recalling that  $ \Omega_{|h|} = \{ x \in \Omega : d(x, \partial \Omega) > h\},$  by the classical difference quotient arguments we obtain $C>0$, independent of $h$, such that
\[
\int_{\Omega_{|h|}}|D^h v|^2  \leq C \int_\Omega | \nabla v|^2, \quad \forall v \in W^{1,2}.\\
\]
Therefore, from \eqref{2lemma1} and \eqref{5lemma1} we conclude that
\begin{equation*}
[[u]]_{\mathcal{N}^{1+\frac{2}{r},r}}^r \leq C \int_\Omega|\nabla u|^{r-2}|D^2 u| ^2,
\end{equation*}
and then we obtain the desired estimate.
\end{proof}

We remark that in addition to Lemma \ref{normalemma}, in order to prove that $u\in \mathcal{N}^{1+\frac{2}{r},r}$, it would remain to guarantee that $u\in W^{1,r}$. However, as an alternative, it is possible to combine the embeddings in Lemma \ref{compactness} with standard Gagliardo--Nirenberg interpolations, assuring that $u$ already belongs to this Nikolskii space under the assumptions of the previous lemma.

Finally, we investigate a classic example of solution to a degenerate P.D.E. In turn, our purpose is to illustrate the connections between the spaces $\mathcal{N}^{\sigma,r}$ with a linearized version of problem \eqref{p1}.
\subsection*{ Example}
\label{exampleapp}
Consider the problem
\begin{equation}
\label{eq1appendix}\tag{$D$}
\begin{cases}
\begin{array}{rll}
-\Delta_p u+u &=& f \mbox{ in } \Omega
\\
\dfrac{\partial u}{\partial \eta} & =&0 \mbox{ on } \partial \Omega.
\end{array}
\end{cases}
\end{equation}

We shall investigate the regularity for solutions of the particular case where $\Omega=B(0,1)$, the unit ball.
Further, in order to avoid additional technical details and for the sake of simplicity, we will work with the case $f\in L^s$, for $s=2$.

Within this framework, on one hand, recall that the regularity results discussed on the text guarantee that any solution of \eqref{eq1appendix} satisfies
\[u\in \mathcal{N}^{1+\frac{2}{r_2},r_2}=\mathcal{N}^{p^\prime,\frac{2p}{p^\prime}},\]
since  for $r_s=s(p-2)+2$ and $s=2$, we have $r_2=\frac{2p}{p^\prime}$.

On the other hand, let us define
\[u(x)=|x|^\alpha-\frac{\alpha}{2}|x|^2,\]
where
\[p^\prime-\frac{N}{2(p-1)}<\alpha \leq 2-\frac{N}{2(p-1)}.\]

For this choice of $u$, we set $f=-\Delta_p u + u$.

We claim that, under the latter assumptions, $u$ is the solution of $\eqref{eq1appendix}$.

Indeed, first remark that
\[\frac{\partial u}{\partial x_i}=\alpha (|x|^{\alpha-2}-1)x_i, \mbox{ for }i=1,\cdots,N, \]
and
\[
\frac{\partial^2 u}{\partial x_j \partial x_i} =
\begin{cases}
\alpha(\alpha-2)|x|^{\alpha-4}x_i x_j \mbox{ if } i\neq j
\\
\alpha (|x|^{\alpha-2}-1)+\alpha(\alpha-2)|x|^{\alpha-4}x_i^2 \mbox{ if } i=j.
\end{cases}
\]
It is then clear that
\[|\nabla u|\leq C(|x|^{\alpha-1}+|x|) \mbox{ and } |D^2 u|\leq C(|x|^{\alpha-2}+1),\]
and thus, after straightforward calculations, we arrive at
\[|\Delta_p u|\leq C(|x|^{(\alpha-2)(p-2)+\alpha+p-4}+|x|^{p-2}+|x|^{\alpha+p-4}).\]
Hence, by combining the choice for $\alpha$ and the Coarea formula, we find out that
\[\int_\Omega |\Delta_p u|^2<+\infty \mbox{ and } \int_\Omega |u|^2 <+\infty.\]
Indeed, it is enough to remark that
\[
\alpha>\max\bigg\{p^\prime -\frac{N}{2(p-1)},4-p-\frac{N}{2},-\frac{N}{2}\bigg\}=p^\prime -\frac{N}{2(p-1)},
\]
since $p>2$.

In this fashion, as $\dfrac{\partial u}{\partial\eta}=0$ on $\partial \Omega$, $u$ is a solution to \eqref{eq1appendix}.

Now, in an analogous manner, by the choice of $\alpha $ we observe that \[u\in W^{1,r_2}=W^{1,\frac{2p}{p^\prime}},\] so that in order to assure that \[u\in \mathcal{N}^{1+\frac{2}{r_2},r_2}=\mathcal{N}^{p^\prime,\frac{2p}{p^\prime}}\]
 it is sufficient to prove
\begin{equation}
\label{eq2appendix}
\int_\Omega |\nabla u|^{r_2-2}|D^2 u|^2<+\infty,
\end{equation}
cf. Lemma \ref{normalemma}.

For this, let us  stress that, as $p>2$ and $\Omega=B(0,1)$, it is obvious that
\[|\nabla u|^{2(p-2)}|D^2u|^2\leq C ( |x|^{2(\alpha-1)(p-2)+2(\alpha-2)}+|x|^{2(\alpha-1)(p-2)}+|x|^{2(\alpha+p-4)}+1)\]

Thence, by recalling that $r_2-2=2(p-2)$, the Coarea formula guarantees the validity of \eqref{eq2appendix}, for
\[
\alpha>\max\bigg\{p^\prime -\frac{N}{2(p-1)},1 -\frac{N}{2(p-1)},4-p-\frac{N}{2}\bigg\}=p^\prime -\frac{N}{2(p-1)}.
\]

On the other hand, direct computations show that
\[|D^2 u|\geq C |x|^{\alpha-2},\]
and then, by the Coarea formula we conclude that \[u \notin W^{2,r_2}, \mbox{ since }\alpha \leq 2-\frac{N}{2(p-1)}.\] Hence the solution indeed does not reach the integer regularity which would generalize the linear case.

Finally, let us remark that, in an analogous manner, if we have asked in addition that $p>\frac{N}{2}$ and considered
\[ p^\prime -\frac{N}{2(p-1)} < \alpha \leq 2-\frac{N}{p}\]
then, we would also have \[u \notin W^{2,p},\] and therefore the second derivatives of the solution do not even reach the ``natural" order integrability $p$, once again contrasting the regularity for the nondegenerate case.

\subsection*{Acknowledgements} The authors would like to express their gratitude to G. M. de Figueiredo 
for his 
fruitful  suggestions to the present work.


\bibliographystyle{amsplain}

\end{document}